\documentclass{article}
\usepackage{theorem,amsmath,amssymb}
\newtheorem{Th}{Theorem}[section]
\newtheorem{Le}{Lemma}[section]
\newtheorem{Co}{Corollary}[section]

\date{}
\title{On the nonstationary Stokes system in a cone}
\author{Vladimir Kozlov and J\"urgen Rossmann}
\begin{document}
\maketitle

\abstract{The authors consider the Dirichlet problem for the nonstationary Stokes system in a threedimensional cone.
They obtain existence and uniqueness results for solutions in weighted Sobolev spaces and prove a regularity
assertion for the solutions. \\


}

\section*{Introduction}

The present paper deals with the Dirichlet problem for the nonstationary Stokes system in a threedimensional cone $K$.
This means, we consider the problem
\begin{eqnarray} \label{stokes1}
&&\frac{\partial u}{\partial t} - \Delta u + \nabla p = f, \quad -\nabla\cdot u = g \ \mbox{ in } K \times (0,\infty), \\ \label{stokes2}
&& u(x,t)=0 \ \mbox{ for }x\in \partial K, \ t>0, \quad u(x,0)=0 \ \mbox{ for }x\in K.
\end{eqnarray}
The goal of the paper is to prove the existence and uniqueness of solutions in weighted Sobolev spaces and
a regularity assertion for the solutions.
A theory for the heat equation in domains with conical points and edges was established in a number of papers in the last
30 years, see \cite{dCN-11}, \cite{k/m-87}, \cite{k/r-11}, \cite{k/r-12}, \cite{nazarov-01}, \cite{Sol-84}, \cite{Sol-01}. This theory
involves in particular existence and uniqueness results for solutions in weighted Sobolev and H\"older spaces, regularity assertions
and the asymptotics of the solutions near vertices and edges. A class of general parabolic problems in a cone was studied
in \cite{kozlov-88}, \cite{kozlov-89}, \cite{kozlov-91}. However, this class of problems does not include
the Stokes system. Although the stationary Stokes system in domains with conical points and edges is well-studied (see, e.~g.,
\cite{Dauge-89}, \cite{mp-83}, \cite{mr-10}), there is still no theory for the nonstationary Stokes system
in domains with singular boundary points. The present paper is a first step in developing such a theory.

An essential part of the paper (Sections 1 and 2) consists of the investigation of the parameter-depending problem
\begin{equation} \label{par1}
s\, \tilde{u} - \Delta \tilde{u} + \nabla \tilde{p} = \tilde{f}, \quad -\nabla\cdot \tilde{u} = \tilde{g}  \
   \mbox{ in } K, \quad \tilde{u} = 0 \ \mbox{ on } \partial K,
\end{equation}
which arises after the Laplace transformation with respect to the time $t$. In Section 1, we prove that this problem
has a uniquely determined variational solution $(\tilde{u},\tilde{p}) \in E_\beta^1(K) \times \big( V_\beta^0(K) + V_\beta^1(K)\big)$ if
$\mbox{Re}\, s \ge 0$, $s\not=0$ and $|\beta|$ is sufficiently small.
Here $V_\beta^l(K)$ denotes the weighted Sobolev space of all functions (vector-functions) with finite norm
\begin{equation} \label{Vbeta}
\| u\|_{V_\beta^l(K)} = \Big( \int_K \sum_{|\alpha|\le l} r^{2(\beta-l+|\alpha|)}\big| \partial_x^\alpha u(x)\big|^2\, dx\Big)^{1/2},
\end{equation}
while $E_\beta^l(K)$ is the  weighted Sobolev space with the norm
\begin{equation} \label{Ebeta}
\| u\|_{E_\beta^l(K)} = \Big( \int_K \sum_{|\alpha|\le l} \big( r^{2\beta}+r^{2(\beta-l+|\alpha|)}\big)\big|
  \partial_x^\alpha u(x)\big|^2\, dx\Big)^{1/2},
\end{equation}
$r=r(x)$ denotes the distance of the point $x$ from the vertex of the cone.

The goal of Section 2 is to prove the existence and uniqueness of solutions $(\tilde{u},\tilde{p})
\in E_\beta^2(K) \times V_\beta^1(K)$ of the parameter-depending problem in the case $\mbox{Re}\, s \ge 0$, $s\not=0$.
Note that the problem (\ref{par1}) is \underline{not} elliptic with parameter in the sense of \cite{av-64}.
Therefore, the results concerning general parabolic problems in a cone which were obtained in
\cite{kozlov-88,kozlov-89,kozlov-91} are not applicable to our problem.

We get two $\beta$-intervals for which we have an existence and uniqueness result in the space $E_\beta^2(K)\times V_\beta^1(K)$,
namely the intervals
\begin{equation} \label{interval}
\frac 12 - \lambda_1^+ <\beta < \frac 12 \quad\mbox{and}\quad  \frac 12 < \beta < \min\Big(\mu_2^+ +\frac 12\, , \lambda_1^+ +\frac 32\Big)
\end{equation}
(see Theorems \ref{t3} and \ref{t4}).
Here, $\lambda_1^+$ and $\mu_2^+$ are positive numbers depending on the cone. More precisely, $\lambda_1^+$ is
the smallest positive eigenvalues of the operator pencil ${\cal L}(\lambda)$ generated by the Dirichlet problem for
the stationary Stokes system, while $\mu_2^+$ is the smallest positive eigenvalue of the operator pencil
${\cal N}(\lambda)$ generated by the Neumann problem for the Laplacian, respectively.
We show that the above inequalities for $\beta$ are sharp.

We also prove a regularity assertion for the solutions of the problem (\ref{par1}), i.~e., we answer the question
under wich conditions a solution $(u,p)\in E_\beta^2(K)\times V_\beta^1(K)$ belongs to the space
$E_\gamma^2(K)\times V_\gamma^1(K)$. In the case $\beta>\gamma$, we get a condition concerning the eigenvalues of
the pencil ${\cal L}(\lambda)$, while the eigenvalues of the pencil ${\cal N}(\lambda)$ must be considered if
$\beta<\gamma$ (see Lemma \ref{l5b}). The appearance of two different operator pencil is one important feature of the solvability
and regularity theory for the problem (\ref{par1}). In analogous results for parameter-elliptic problems
(see \cite{kozlov-88,kozlov-89}),  only the eigenvalues of one operator pencil play a role.
A further important difference with the parameter-elliptic case is the use of different function spaces.
In the parameter-elliptic case, all functions (solutions and right-hand sides) belong to the spaces $E_\beta^l$.
When considering the problem (\ref{par1}), we seek solutions $(\tilde{u},\tilde{p}) \in E_\beta^2(K)\times V_\beta^1(K)$,
while the right-hand sides $\tilde{f}$, $\tilde{g}$ belong to the spaces $E_\beta^0(K)$ and
$X_\beta^1(K)=V_\beta^1(K)\cap \big(V_{-\beta}^1(K)\big)^*$, respectively.

In Section 3, we apply the results concerning the parameter-depending problem and obtain the following
existence and uniqueness result for solutions of the problem (\ref{stokes1}), (\ref{stokes2}):
If the data $f$ and $g$ belong to corresponding weighted Sobolev spaces, $\beta$ lies in one of the intervals
(\ref{interval}) and $\int_K g(x,t)\, dx=0$ for almost all $t$ if $\beta>1/2$,
then there exists a unique solution $(u,p)$ of the problem (\ref{stokes1}), (\ref{stokes2}) such that
\[
u \in L_2\big({\Bbb R}_+,V_\beta^2(K)\big), \ \ \partial_t u \in L_2\big({\Bbb R}_+,V_\beta^0(K)\big), \quad
p \in L_2\big({\Bbb R}_+,V_\beta^1(K)\big).
\]
Moreover, a regularity assertion for this solution is proved.

\section{Weak solutions of the parameter-dependent problem}

Let $s$ be an arbitrary complex number, $\mbox{Re}\, s \ge 0$. In this section, we prove the existence and
uniqueness of weak solutions of the boundary value problem
\begin{equation} \label{par3}
s\, u - \Delta u + \nabla p = f, \ \ - \nabla\cdot u = g \ \mbox{ in }K, \quad u=0\ \mbox{ on }\partial K\backslash \{ 0\}
\end{equation}
in weighted Sobolev spaces.

\subsection{Weighted Sobolev spaces on the cone}

Let $K = \big\{ x\in {\Bbb R}^3:\ \omega=x/|x| \in \Omega\}$ be a cone with vertex at the origin. Here,
$\Omega$ is a subdomain of the unit sphere $S^2$ with smooth (of class $C^2$) boundary $\partial\Omega$.
For nonnegative integer $l$ and real $\beta$, we define the weighted Sobolev spaces $V_\beta^l(K)$ and $E_\beta^l(K)$
as the sets of all functions (or vector functions) with finite norms (\ref{Vbeta}) and (\ref{Ebeta}),
respectively.
Note that the spaces $V_\beta^l(K)$ and $E_\beta^l(K)$ can be also defined as the closures of
$C_0^\infty(\overline{K}\backslash \{ 0\})$ with respect to the above norms. The space $E_\beta^l(K)$
can be also equipped with the equivalent norm
\[
\| u\| = \Big( \int_K r^{2\beta} \sum_{|\alpha|=l} \big| \partial_x^\alpha u(x)\big|^p\, dx
  + \int_K \big( r^{2\beta}+r^{2(\beta-l)}\big)\, \big| u(x)\big|^p\, dx\Big)^{1/2}.
\]
In the case $l\ge 1$, the trace spaces for $V_\beta^l(K)$ and $E_\beta^l(K)$ on
$\partial K \backslash\{ 0\}$ are denoted by $V_\beta^{l-1/2}(\partial K)$ and $E_\beta^{l-1/2}(\partial K)$,
respectively. If $l\le 0$, then $V_\beta^{l-1/2}(\partial K)$ denotes the dual space of
$V_{-\beta}^{-l+1/2}(\partial K)$.
Furthermore, we define $\stackrel{\circ}{V}\!{}_\beta^1(K)$ and $\stackrel{\circ}{E}\!{}_\beta^1(K)$
as the spaces of all functions $u\in V_\beta^1(K)$ and $u\in E_\beta^1(K)$, respectively, which are zero on
$\partial K\backslash \{ 0\}$. The dual spaces of $\stackrel{\circ}{V}\!{}_\beta^1(K)$ and $\stackrel{\circ}{E}\!{}_\beta^1(K)$
are denoted by $V_{-\beta}^{-1}(K)$ and $E_{-\beta}^{-1}(K)$, respectively.
Since
\[
\int_K r^{2\beta-2} \, \big| u(x)\big|^2\, dx \le c\, \int_K r^{2\beta}\, \big| \nabla u(x)\big|^2\, dx
\]
for $u\in C_0^\infty(K)$, the norm
\[
\| u \| = \Big( \int_K r^{2\beta}\, \big( |u|^2 + |\nabla u|^2\big)\, dx\Big)^{1/2}
\]
is equivalent to the $E_\beta^1(K)$-norm in $\stackrel{\circ}{E}\!{}_\beta^1(K)$. Finally, we mention that there
exists a constant $c$ such that
\begin{equation} \label{estV}
\int_{\partial K} r^{2(\beta-l)+1}\, \big| u(x)\big|^2\, dx \le c\, \| u\|^2_{V_\beta^{l-1/2}(\partial K)}
\end{equation}
for all $u \in V_\beta^{l-1/2}(\partial K)$, $l\ge 1$, and
\begin{equation} \label{estE}
\int_{\partial K} \big( r^{2\beta} + r^{2(\beta-l)+1}\big)\, \big| u(x)\big|^2\, dx \le c\, \| u\|^2_{E_\beta^{l-1/2}(\partial K)}
\end{equation}
for all $u \in E_\beta^{l-1/2}(\partial K)$, $l\ge 1$ (cf. \cite[Lemmas 1.4 and 1.5]{mp-78}).

\subsection{A property of the operator div}

It is obvious that the operator div realizes a linear and continuous mapping
\[
\stackrel{\circ}{E}\!{}_0^1(K) \to L_2(K) \cap \big( V_0^1(K)\big)^*
\]
Our goal is to show that this operator is surjective. To this end, we introduce the following operator pencil
${\cal L}$ generated by the stationary Stokes system in the cone $K$. For every complex $\lambda$,
we define the operator ${\cal L}(\lambda)$ as the mapping
\begin{eqnarray*}
&& \stackrel{\circ}{W}\!{}^1(\Omega)\times L_2(\Omega) \ni \left( \begin{array}{c} u \\ p \end{array}\right) \\ && \qquad
  \to \left( \begin{array}{c} r^{2-\lambda}\big(-\Delta r^{\lambda}u(\omega)+\nabla r^{\lambda-1}p(\omega)\big)\\*[1ex]
  -r^{1-\lambda}\nabla\cdot\big( r^\lambda u(\omega)\big)\end{array}\right) \in W^{-1}(\Omega)\times L_2(\Omega),
\end{eqnarray*}
where $r=|x|$ and $\omega=x/|x|$. The properties of the pencil ${\cal L}$ are studied, e.g., in \cite{kmr-01}.
In particular, it is known that the eigenvalues of this pencil in the strip $-2\le \mbox{Re}\, \lambda \le 1$
are real and that the numbers $1$ and $-2$ are always eigenvalues. If $\Omega$ is contained in a half-sphere,
then these numbers are the only eigenvalues in the interval $[-2,1]$ (cf. \cite[Theorem 5.5.5]{kmr-01}).
We use this fact for the proof of the following lemma.

\begin{Le} \label{l1}
For arbitrary $g\in L_2(K)\cap (V_0^1(K))^*$, there exists a vector function $v\in \stackrel{\circ}{E}\!{}_0^1(K)$
sucht that \ $-\nabla\cdot v =g$ in $K$ and
\begin{equation} \label{1l1}
\| v\|_{E_0^1(K)} \le c\, \Big( \| g\|_{L_2(K)} + \| g\|_{(V_0^1(K))^*} \Big).
\end{equation}
Here, $c$ is a constant independent of $g$.
\end{Le}

P r o o f.
Assume first that $\Omega$ is contained in a half-sphere. We denote the operator $(u,p) \to (f,g)$ of the boundary value problem
\begin{equation} \label{statstokes}
- \Delta u + \nabla p = f, \ \ - \nabla\cdot u = g \ \mbox{ in }K, \quad u=0\ \mbox{ on }\partial K\backslash \{ 0\}
\end{equation}
by $A_0$. Furthermore, we introduce the spaces
\begin{eqnarray*}
&& X_1 = \stackrel{\circ}{V}\!{}_0^1(K) \times L_2(K), \quad Y_1 = V_0^{-1}(K) \times L_2(K), \\
&& X_2 = \big( V_0^2(K) \cap \stackrel{\circ}{V}\!{}_{-1}^1(K)\big) \times V_0^1(K), \quad Y_2 = L_2(K) \times V_0^1(K).
\end{eqnarray*}
By \cite[Theorem 10.2.11]{mr-10}, the operator $A_0$ realizes an isomorphism $X_1 \to Y_1$, i.e., for every
pair $(f,g) \in Y_1$, there exists a uniquely determined vector function $(u,p) \in X_1$ such that
\[
b_0(u,v) - \int_K p\, \nabla\cdot v\, dx = \langle f,v \rangle \ \mbox{ for all }v\in \stackrel{\circ}{V}\!{}_0^1(K), \quad
  -\nabla\cdot u = g \ \mbox{ in }K,
\]
where $b_0(u,v)$ is the bilinear form
\[
b_0(u,v) = \sum_{j=1}^3 \int_K \nabla u_j \cdot \nabla v_j \, dx
\]
(here $u_j$, $v_j$ are the components of the vector functions $u$ and $v$, respectively).
Furthermore by \cite[Theorem 10.3.2]{mr-10}, the operator $A_0$ realizes an isomorphism $X_2 \to Y_2$. From
\cite[Theorem 10.6.9]{mr-10} it follows that every solution $(u,p) \in X_1$ of the equation
\[
A_0\, (u,p) = (f,g)
\]
belongs to the space $X_2$ if $(f,g) \in Y_1\cap Y_2$. Consequently, the operator $A_0$ is also an isomorphism
from $X_1\cap X_2$ onto $Y_1\cap Y_2$ and from $X_1+X_2$ onto $Y_1+Y_2$,
and the adjoint operator $A_0^*$ realizes an isomorphism $Y_1^* \cap Y_2^* \to X_1^* \cap X_2^*$.
This means that for every $(f,g) \in X_1^* \cap X_2^*$, there exists a uniquely determined element
$(v,q) \in Y_1^* \cap Y_2^*$ such that
\begin{eqnarray*}
&& (u,f)_K + (p,g)_K = (-\Delta u + \nabla p, v)_K - (\nabla\cdot u, q)_K \ \mbox{ for all }(u,p) \in X_2, \\
&& (u,f)_K + (p,g)_K = b_0(u,v) - (p, \nabla\cdot v)_K - (\nabla\cdot u, q)_K \ \mbox{ for all }(u,p) \in X_1.
\end{eqnarray*}
Here $(\cdot,\cdot)_K$ denotes the extension of the $L_2(K)$ scalar product to pairs from $X\times X^*$.
Moreover, there exists a constant $c$ independent of $f,g$ such that
\[
\| (v,q)\|_{Y_1^* \cap Y_2^*} \le c\, \| (f,g) \|_{X_1^* \cap X_2^*} \, .
\]
In particular (for $f=0$, $u=0$), we conclude that for arbitrary $g \in L_2(K) \cap V_0^1(K)^*$, there exists
a vector function $v\in \stackrel{\circ}{V}\!{}_0^1(K) \cap L_2(K)=\stackrel{\circ}{E}\!{}_0^1(K)$ such that
\[
(p,g)_K = (\nabla p,v)_K \ \mbox{ for all }p\in V_0^1(K)\quad \mbox{and}\quad
(p,g)_K = -(p,\nabla\cdot v)_K  \ \mbox{ for all }p\in L_2(K).
\]
Furthermore, $v$ satisfies (\ref{1l1}). This proves the assertion for the case when $\Omega$ is contained in a
half-sphere.

In the contrary case, one can find subdomains $\Omega_j\subset\Omega$, $j=1,\ldots,n$, with smooth boundaries and
smooth function $\chi_j$ on $\Omega$ such that $\Omega=\Omega_1\cup\cdots\cup \Omega_n$, $\Omega_j$ is contained
in a half-sphere for every $j$,
\[
\mbox{supp}\, \chi_j \subset\overline{\Omega}_j\quad\mbox{and}\quad \sum_{j=1}^n \zeta_j =1.
\]
Let $K_j$ be the cone consisting of all $x\in K$ such that $x/|x| \in \Omega_j$.
By the first part of the proof, there exist vector functions $w_j \in \stackrel{\circ}{E}\!{}_0^1(K_j)$
such that $\nabla \cdot w_j(x) = -\chi_j(\omega)\, g(x)$ for $x\in K_j$ and
\[
\| w_j\|_{E_0^1(K_j)} \le c\, \Big( \| \chi_j g\|_{L_2(K_j)} + \| \chi_j g\|_{(V_0^1(K_j))^*} \Big)
\]
Let $v_j$ be the extensions of $w_j$ by zero to the cone $K$. It is obvious that $v_j \in \stackrel{\circ}{E}\!{}_0^1(K)$
and that the sum $v=v_1+\cdots + v_m$ satisfies the equation $\nabla\cdot v=-g$ in $K$ and the estimate (\ref{1l1}).
\hfill $\Box$

\subsection{Existence of variational solutions}

In order to define variational solutions of the problem (\ref{par3}), we introduce the
bilinear form
\[
b_s(u,v) = \int_K \Big( su\cdot v + \sum_{j=1}^3 \nabla u_j \cdot \nabla v_j\Big)\, dx.
\]
It is clear that the mappings
\[
v \to b_s(u,v) \quad \mbox{and} \quad v \to \int_K p\, \nabla\cdot v\, dx
\]
define linear and continuous functionals on $\stackrel{\circ}{E}\!{}_0^1(K)$ if $u\in \stackrel{\circ}{E}\!{}_0^1(K)$
and $p\in L_2(K)+V_0^1(K)$. By a {\em variational solution} of the problem (\ref{par3}), we mean a pair
$(u,p) \in \stackrel{\circ}{E}\!{}_0^1(K) \times \big( L_2(K)+V_0^1(K)\big)$ satisfying the equations
\begin{equation} \label{1t1}
b_s(u,v) - \int_K p\, \nabla\cdot v\, dx = \langle f,v \rangle \ \mbox{ for all }v\in \stackrel{\circ}{E}\!{}_0^1(K), \qquad
-\nabla\cdot u = g \ \mbox{ in }K,
\end{equation}
where $f\in E_0^{-1}(K)$ and $g\in L_2(K)\cap \big(V_0^1(K)\big)^*$ are given.

\begin{Th} \label{t1}
Suppose that $|s|=1$, $\mbox{\em Re}\, s \ge 0$, $f\in E_0^{-1}(K)$ and $g\in L_2(K)\cap \big(V_0^1(K)\big)^*$.
Then there exists a uniquely determined solution $(u,p) \in \stackrel{\circ}{E}\!{}_0^1(K) \times \big( L_2(K)+V_0^1(K)\big)$
of the problem {\em (\ref{1t1})}. Furthermore,
\[
\| u\|_{E_0^1(K)} + \| p\|_{L_2(K)+V_0^1(K)} \le c\, \Big( \| f\|_{E_0^{-1}(K)} + \| g\|_{L_2(K)} + \| g\|_{(V_0^1(K))^*}\Big)
\]
with a constant $c$ independent of $f,g$ and $s$.
\end{Th}

P r o o f.
First, we prove the existence of solutions in the case $g=0$. We introduce the subspace
\[
H = \{ u\in \stackrel{\circ}{E}\!{}_0^1(K): \ \nabla\cdot u =0 \ \mbox{in }K\}
\]
of $\stackrel{\circ}{E}\!{}_0^1(K)$ and denote by $H^\bot$ its orthogonal complement in $\stackrel{\circ}{E}\!{}_0^1(K)$.
Since
\[
\big| b_s(u,\bar{u})\big| \ge \frac 12 \, \| u\|^2_{E_0^1(K)} \ \mbox{ for all }u\in E_0^1(K)
\]
if $|s|=1$, $\mbox{Re}\, s\ge 0$, there exists a uniquely determined vector function $u\in H$ such that
\begin{equation} \label{2t1}
b_s(u,v)= \langle f,v \rangle\quad \mbox{for all }v\in H, \qquad \| u\|_{E_0^1(K)} \le c\, \| f\|_{E_0^{-1}(K)}\, .
\end{equation}
By Lemma \ref{l1}, the operator $\mbox{div}$ is an isomorphism from $H^\bot$ onto $L_2(K)\cap \big(V_0^1(K)\big)^*$,
i.e., for every $q \in L_2(K)\cap \big(V_0^1(K)\big)^*$, there exists a uniquely determined vector function
$v=(-\mbox{div})^{-1}q \in H^\bot$ such that $\nabla\cdot v=-q$ and
\[
\| v\|_{E_0^1(K)} \le c\, \Big( \| q\|_{L_2(K)} + \| q\|_{(V_0^1(K))^*}\Big) .
\]
We consider the functional
\begin{equation} \label{3t1}
\ell(q) =  \langle f,v \rangle- b_s(u,v) = \langle f, (-\mbox{div})^{-1}q\rangle - b_s\big(u,(-\mbox{div})^{-1}q\big)
\end{equation}
on $L_2(K)\cap \big(V_0^1(K)\big)^*$. Obviously,
\begin{eqnarray*}
\big| \ell(q)\big| & \le & c\, \Big( \| f\|_{E_0^{-1}(K)} + \| u\|_{E_0^1(K)}\Big)\, \| v\|_{E_0^1(K)} \\
& \le & c\, \| f\|_{E_0^{-1}(K)}\, \Big( \| q\|_{L_2(K)} + \| q\|_{(V_0^1(K))^*}\Big).
\end{eqnarray*}
This means that $\ell$ is continuous on $L_2(K)\cap \big(V_0^1(K)\big)^*$ and there exists an element $p$ of the dual
space $L_2(K)+V_0^1(K)$ such that $\ell(q) =  \int_K p\, q dx$ for all $q\in L_2(K)\cap \big(V_0^1(K)\big)^*$.
Then by (\ref{3t1}), we have
\[
b_s(u,v) + \int_K p\, q\, dx = \langle f,v \rangle
\]
for $q\in L_2(K)\cap \big(V_0^1(K)\big)^*$, $v=(-\mbox{div})^{-1}q \in H^\bot$. This means that
\[
b_s(u,v) -  \int_K p \, \nabla \cdot v\, dx =  \langle f,v \rangle \ \mbox{ for all }v\in H^\bot.
\]
Using (\ref{2t1}), we conclude that $(u,p)$ is a solution of the problem (\ref{1t1}) in the case $g=0$.

Now let $g\in L_2(K)\cap \big(V_0^1(K)\big)^*$ be an arbitrary function. By Lemma \ref{l1}, there exists
a function $w\in \stackrel{\circ}{E}\!{}_0^1(K)$ such that
\[
- \nabla\cdot w = g \, \mbox{ in }K, \quad \| w\|_{E_0^1(K)} \le c\, \Big( \| g \|_{L_2(K)} + \| g\|_{(V_0^1(K))^*}\Big) .
\]
Obviously the mapping
\[
\stackrel{\circ}{E}\!{}_0^1(K) \ni v \to \langle f,v \rangle - b_s(w,v)
\]
defines a continuous functional on $\stackrel{\circ}{E}\!{}_0^1(K)$. Therefore, by the first part of the proof,
there exist functions $W\in \stackrel{\circ}{E}\!{}_0^1(K)$ and $p\in L_2(K)+V_0^1(K)$ such that
\[
b_s(W,v) - \int_K p\, \nabla\cdot v \, dx = \langle f,v \rangle - b_s(w,v) \  \mbox{ for all }v\in \stackrel{\circ}{E}\!{}_0^1(K).
\]
Then $(u,p)=(W+w,p)$ is a solution of the problem (\ref{1t1}).

We prove the uniqueness of the solution. Let $(u,p)\in \stackrel{\circ}{E}\!{}_0^1(K)\times
\big( L_2(K)+V_0^1(K)\big)$ be a solution of the problem (\ref{1t1}) with $f=0$ and $g=0$. Then in particular,
$u\in H$ and $b_s(u,\bar{u})=0$ what implies $u=0$ and
$\int_K p\, \nabla\cdot v\, dx =0$ for all $v\in \stackrel{\circ}{E}\!{}_0^1(K)$. Using Lemma \ref{l1}, we conclude
that $(p,q)_K=0$ for all $q\in L_2(K)\cap \big(V_0^1(K)\big)^*$ and, consequently, $p=0$.
The proof of the theorem is complete. \hfill $\Box$ \\

As a consequence of the last theorem, we obtain the existence and uniqueness of variational solutions in the space
$\stackrel{\circ}{E}\!{}_\beta^1(K) \times \big( V_\beta^0(K)+V_\beta^1(K)\big)$ if $|\beta|$ is
sufficiently small. Suppose that $f\in E_\beta^{-1}(K)$ and $g\in V_\beta^0(K)\cap (V_{-\beta}^1(K))^*$.
Then $(u,p) \in \stackrel{\circ}{E}\!{}_\beta^1(K) \times \big( V_\beta^0(K)+V_\beta^1(K)\big)$
is called a variational solution of the problem (\ref{par3}) if
\begin{equation} \label{varsol}
b_s(u,v) - \int_K p\, \nabla\cdot v\, dx = \langle f,v \rangle \ \mbox{ for all }v\in \stackrel{\circ}{E}\!{}_{-\beta}^1(K), \qquad
-\nabla\cdot u = g \ \mbox{ in }K.
\end{equation}

\begin{Co}
Suppose that $|s|=1$, $\mbox{\em Re}\, s \ge 0$, $f\in E_\beta^{-1}(K)$ and $g\in V_\beta^0(K)\cap \big(V_{-\beta}^1(K)\big)^*$,
where $|\beta|$ is sufficiently small.
Then there exists a uniquely determined solution $(u,p) \in \stackrel{\circ}{E}\!{}_\beta^1(K) \times \big( V_\beta^0(K)+V_\beta^1(K)\big)$
of the problem {\em (\ref{varsol})}. This solution satisfies the estimate
\[
\| u\|_{E_\beta^1(K)} + \| p\|_{V_\beta^0(K)+V_\beta^1(K)} \le c\, \Big( \| f\|_{E_\beta^{-1}(K)}
  + \| g\|_{V_\beta^0(K)} + \| g\|_{(V_{-\beta}^1(K))^*}\Big)
\]
with a constant $c$ independent of $f,g$ and $s$.
\end{Co}

P r o o f.
Let $(U,P) = r^\beta\, (u,p)$. Obviously,
\[
(u,p) \in \stackrel{\circ}{E}\!{}_\beta^1(K) \times \big( V_\beta^0(K)+V_\beta^1(K)\big)\ \Leftrightarrow \
(U,P) \in \stackrel{\circ}{E}\!{}_0^1(K) \times \big( L_2(K)+V_0^1(K)\big).
\]
Furthermore, $(u,p)$ is a solution of (\ref{varsol}) if and only if $(U,P)$ satisfies
\begin{eqnarray} \label{varsol1}
&& b_s(r^{-\beta}U,r^{\beta}V) - \int_K r^{-\beta} P\, \nabla\cdot r^\beta V\, dx = \langle F,V \rangle\
  \mbox{ for all }V=r^{-\beta}v\in \stackrel{\circ}{E}\!{}_0^1(K), \\ \label{varsol2}
&& -r^\beta \nabla\cdot (r^{-\beta} U) = G \ \mbox{ in }K,
\end{eqnarray}
where $F=r^\beta f \in E_0^{-1}(K)$ and $G=r^\beta g\in L_2(K)\cap \big(V_0^1(K)\big)^*$. We denote the operator
\[
\stackrel{\circ}{E}\!{}_0^1(K) \times \big( L_2(K)+V_0^1(K)\big) \ni (U,P) \to (F,G) \in E_0^{-1}(K)\times \big( L_2(K)\cap \big(V_0^1(K)\big)^*\big)
\]
by $A_\beta$. By Theorem \ref{t1}, the operator $A_0$ is an isomorphism. Since the operator $A_0-A_\beta$ has a small norm
for small $|\beta|$, it follows that $A_\beta$ is an isomorphism if $|\beta|$ is sufficiently small. This proves the
corollary. \hfill $\Box$

\section{Strong solutions of the problem (\ref{par3})}

Now, we are interested in solutions $(u,p) \in E_\beta^2(K) \times V_\beta^1(K)$ of the problem (\ref{par3}).

\subsection{The operator of the problem (\ref{par3})}

Obviously, the mapping
\[
E_\beta^2(K) \times V_\beta^1(K) \ni (u,p) \to f = su-\Delta u + \nabla p \in E_\beta^0(K)
\]
and the mapping $E_\beta^2(K) \ni u\to g=-\nabla\cdot u\in E_\beta^1(K)$ are continuous for arbitrary real $\beta$.
Furthermore, for arbitrary  $u\in E_\beta^2(K) \cap \stackrel{\circ}{E}\!{}_\beta^1(K)$
and $q\in V_{-\beta}^1(K)$, we get
\[
\Big| \int_K q\, \nabla\cdot u\, dx \Big| = \Big| \int_K u\cdot \nabla q\, dx\Big| \le \| u\|_{V_\beta^0(K)}\
  \| q\|_{V_{-\beta}^1(K)} \le \| u\|_{E_\beta^2(K)}\ \| q\|_{V_{-\beta}^1(K)} \, .
\]
Thus, the mapping $u\to g=-\nabla\cdot u$ is continuous from $E_\beta^2(K) \cap \stackrel{\circ}{E}\!{}_\beta^1(K)$
into the space
\[
X_\beta^1(K)  =  E_\beta^1(K) \cap \big(V_{-\beta}^1(K)\big)^* .
\]
and the operator of the problem (\ref{par3}) realizes a linear and continuous mapping
from $\big( E_\beta^2(K) \cap \stackrel{\circ}{E}\!{}_{\beta}^1(K)\big) \times V_\beta^1(K)$ into $E_\beta^0(K)
\times X_\beta^1(K)$ for arbitrary real $\beta$. We denote the operator
\begin{eqnarray} \label{1t2}
\big( E_\beta^2(K) \cap \stackrel{\circ}{E}\!{}_{\beta}^1(K)\big) \times V_\beta^1(K) \ni (u,p)  \to \
  \big( su-\Delta u + \nabla p, - \nabla\cdot u\big) \in E_\beta^0(K)\times X_\beta^1(K)
\end{eqnarray}
of this problem by $A_\beta$. Note that
\[
\int_K  r^{2\beta}\, |g|^2\, dx \le \| g\|_{(V_{-\beta}^1(K))^*} \
  \| r^{2\beta} g\|_{V_{-\beta}^1(K)} \le c\, \| g\|_{(V_{-\beta}^1(K))^*} \ \| g\|_{V_\beta^1(K)}\, .
\]
for $g\in V_\beta^1(K)\cap \big( V_{-\beta}^1(K)\big)^*$. Consequently, $V_\beta^0(K) \subset V_\beta^1(K) \cap
\big( V_{-\beta}^1(K)\big)^*$. This implies that
\[
X_\beta^1(K) = V_\beta^1(K) \cap (V_{-\beta}^1(K))^*.
\]

\begin{Le} \label{l1b}
Suppose that $u\in E_\beta^2(K) \cap \stackrel{\circ}{E}\!{}_{\beta}^1(K)$, $p\in
V_\beta^1(K)$ and $A_\beta(u,p)=(f,g)$. If $0\le \beta\le 1$, then $(u,p)$ coincides with the
variational solution of the problem  {\em (\ref{par3})} in the space $\stackrel{\circ}{E}\!{}_0^1(K) \times
\big( L_2(K)+V_0^1(K)\big)$.
\end{Le}

P r o o f.
Suppose that $0\le \beta\le 1$. Then
\[
E_\beta^2(K) \cap \stackrel{\circ}{E}\!{}_{\beta}^1(K) \subset \stackrel{\circ}{E}\!{}_0^1(K) \quad \mbox{ and }
V_\beta^1(K) \subset L_2(K)+V_0^1(K).
\]
Furthermore, $E_\beta^0(K)\subset E_0^{-1}(K)$ and $X_\beta^1(K) \subset L_2(K) \cap \big( V_0^1(K)\big)^*$ for
$0\le \beta \le 1$. Obviously, every solution $(u,p)\in E_\beta^2(K) \times V_\beta^1(K)$ of (\ref{par3})
satisfies (\ref{1t1}). This proves the lemma. \hfill $\Box$ \\

In particular, it follows from Theorem \ref{t1} and Lemma \ref{l1b} that the kernel of $A_\beta$ is trivial for
$0\le \beta \le 1$. Later, we will improve this result.

\subsection{An a priori estimate for the solutions of the problem (\ref{par3})}

In order to prove a local estimate for the solutions of the problem (\ref{par3}), we employ the following
lemma.

\begin{Le} \label{l2}
Let ${\cal D}$ be a bounded domain with smooth (of class $C^2$) boundary $\partial{\cal D}$.
Suppose that $u\in W_2^2({\cal D})$, $\nabla p \in L_2({\cal D})$,
\[
s\, u - \Delta u + \nabla p = f, \ \ -\nabla\cdot u=g \ \mbox{ in }{\cal D}, \quad u=0\ \mbox{ on }\partial{\cal D}.
\]
Then
\begin{equation} \label{1l2}
\| D^2 u\|_{L_2({\cal D})} + |s|\, \| u\|_{L_2({\cal D})} + \| \nabla p\|_{L_2({\cal D})} \le c\, \Big( \| f\|_{L_2({\cal D})}
  + \| g\|_{W_2^1({\cal D})} + |s|\, \| g\|_{W_2^{-1}({\cal D})}+ \| u\|_{L_2({\cal D})}\Big)
\end{equation}
with a constant independent of $s$. Here, $D^2 u$ denotes the vector of all second order derivatives of $u$.
\end{Le}

P r o o f.
The functions $v(x,t)=e^{st} u(x)$, $q(x,t)=e^{st}p(x)$ satisfy the equations
\begin{eqnarray*}
&&\partial_t v(x,t) - \Delta v(x,t) + \nabla q(x,t) = e^{st} f(x), \ \ - \nabla\cdot v(x,t) = e^{st} g(x) \ \mbox{ for }x\in {\cal D},\\
&& v(x,t)=0 \ \mbox{ for }x\in \partial{\cal D}, \quad v(x,0)=u(x) \ \mbox{ for }x\in {\cal D}.
\end{eqnarray*}
Consequently by \cite[Theorem 3.1]{sol-03},
\begin{eqnarray*}
&& \| D^2 v\|_{L_2({\cal D}\times (0,1))} + \|\partial_t v\|_{L_2({\cal D}\times (0,1))} + \| \nabla q\|_{L_2({\cal D}\times (0,1))}\\
&&  \le c\, \Big( \| e^{st}f\|_{L_2({\cal D}\times (0,1))} + \| e^{st}g\|_{L_2(0,1;W_2^1({\cal D}))} +
  \| \partial_t(e^{st}g)\|_{L_2(0,1;W_2^{-1}({\cal D}))} + \| u\|_{L_2({\cal D})}\Big).
\end{eqnarray*}
The result follows. \hfill $\Box$ \\

In the following, let $\zeta_\nu$ be infinitely differentiable functions on ${\Bbb R}_+=(0,\infty)$ with support
in the interval $(2^{\nu-1},2^{\nu+1})$ such that
\[
\big| \partial_r^j \zeta_\nu(r)\big| \le 2^{-j\nu} \ \mbox{ for }j=0,1,\ldots \ \mbox{ and } \ \sum_{\nu=-\infty}^{+\infty} \zeta_\nu=1.
\]
Using the notation $r=|x|$, the functions $x\to \zeta_\nu(r)$ can be considered as functions in $K$.

\begin{Le} \label{l3}
Suppose that $(u,p)\in \big( E_\beta^2(K) \cap \stackrel{\circ}{E}\!{}_{\beta}^1(K)\big) \times V_\beta^1(K)$ is a solution
of {\em (\ref{par3})}. Then
\begin{eqnarray} \label{1l3}
&& \hspace{-2em}\| D^2 (\zeta_\nu u)\|^2_{L_2(K)} + |s|^2\, \| \zeta_\nu u\|^2_{L_2(K)} + \| \nabla (\zeta_\nu p)\|^2_{L_2(K)} \nonumber\\
&& \le c\, \Big( \| f\|^2_{L_2(K_\nu)} + \| \nabla (\zeta_\nu g)\|^2_{L_2(K_\nu)} + 2^{-2\nu}\, \|\zeta_\nu g\|^2_{L_2(K_\nu)}
  + |s|^2\, \| \zeta_\nu g\|^2_{W_2^{-1}(K_\nu)}  \nonumber \\
&& \hspace{6em} + \ 2^{-2\nu}\, \| \nabla u \|^2_{L_2(K_\nu)} + 2^{-4\nu} \, \| u\|^2_{L_2(K_\nu)} + 2^{-2\nu}\, \| p\|^2_{L_2(K_\nu)}\Big),
\end{eqnarray}
where $K_\nu= \{ x\in K:\, 2^{\nu-1}< |x| < 2^{\nu+1}\}.$ The constant $c$ in {\em (\ref{1l3})} is independent of $\nu$ and $s$.
\end{Le}

P r o o f.
Obviously,
\begin{equation} \label{2l3}
s\, \zeta_\nu u - \Delta(\zeta_\nu u) + \nabla(\zeta_\nu p) = F, \ \ -\nabla\cdot(\zeta_\nu u) = G \ \mbox{ in }K,
\end{equation}
where $F_i=\zeta_\nu f_i - 2\nabla\zeta_\nu \cdot \nabla u_i - u_i\, \Delta\zeta_\nu + p\,\partial_{x_i}\zeta_\nu$ for $i=1,2,3$
and $G=\zeta_\nu g - u\cdot \nabla\zeta_\nu$. We define
\[
v(x) = \zeta_\nu(2^\nu r)\, u(2^\nu x),\quad q(x) = 2^\nu \zeta_\nu(2^\nu r)\, p(2^\nu x).
\]
By (\ref{2l3}), we have
\[
2^{2\nu}s\, v - \Delta v + \nabla q = \Phi, \ \ -\nabla\cdot v = \Psi \ \mbox{ in }K,
\]
where $\Phi(x)=2^{2\nu}\, F(2^\nu x)$ and $\Psi(x)=2^\nu\, G(2^\nu x)$. Since $v$ and $q$ vanish
outside the set $K_0=\{ x\in K: 1/2 < |x| < 2\}$, it follows from Lemma \ref{l2} that
\begin{eqnarray} \label{3l3}
&& \| D^2 v\|_{L_2(K)} + 2^{2\nu}\, |s|\, \| v\|_{L_2(K)} + \| \nabla q\|_{L_2(K)}  \nonumber \\
&& \le c\, \Big( \| \Phi\|_{L_2(K_0)} + \| \Psi\|_{W_2^1(K_0)}
  + 2^{2\nu}\, |s|\, \| \Psi\|_{W_2^{-1}(K_0)}+ \| v\|_{L_2(K_0)}\Big),
\end{eqnarray}
where the constant $c$ is independent of $\nu$ and $s$. One can easily check that
\[
\| v\|_{L_2(K)} = 2^{-\nu n/2} \| \zeta_\nu u\|_{L_2(K)}\, , \quad
\| D^2 v\|_{L_2(K)} = 2^{\nu(2-n/2)}\, \| D^2(\zeta_\nu u)\|_{L_2(K)}
\]
and
\[
\| \nabla q\|_{L_2(K)} = 2^{\nu(2-n/2)}\, \| \nabla(\zeta_\nu p\|_{L_2(K)}\, .
\]
Furthermore, we obtain the estimates
\begin{eqnarray*}
&& \hspace{-2em}\| \Phi \|_{L_2(K)} +\| \Psi\|_{W_2^1(K)} = 2^{\nu(2-n/2)}\, \Big( \| F \|_{L_2(K)} + \| \nabla G\|_{L_2(K)}
  + 2^{-\nu}\, \| G\|_{L_2(K)}\Big) \\
&& \le  c\, 2^{\nu(2-n/2)}\, \Big( \| \zeta_\nu f\|_{L_2(K)} + \| \nabla(\zeta_\nu g)\|_{L_2(K)} + 2^{-\nu}\, \| \zeta_\nu g\|_{L_2(K)} \\
&& \hspace{6em}  + \ 2^{-\nu}\, \| \nabla u\|_{L_2(K_\nu)} + 2^{-2\nu}\, \| u\|_{L_2(K_\nu)}
    + 2^{-\nu}\, \|p\|_{L_2(K_\nu)}\Big).
\end{eqnarray*}
Moreover,
\[
\| \Psi\|_{W_2^{-1}(K_0)} = 2^{-\nu n/2} \| G\|_{W_2^{-1}(K_\nu)}
  \le c\, 2^{-\nu n/2} \Big( \| \zeta_\nu g\|_{W_2^{-1}(K_\nu)} + 2^{-\nu} \| u\|_{W_2^{-1}(K_\nu)}\Big)
\]
Using the equality $su=f+\Delta u - \nabla p$ and the estimate $\| f\|_{W_2^{-1}(K_\nu)} \le  2^\nu\, \| f\|_{L_2(K_\nu)}$, we get
\[
\| \Psi\|_{W_2^{-1}(K_0)}  \le c\, 2^{-\nu n/2}  \| \zeta_\nu g\|_{W_2^{-1}(K_\nu)}
  + c\,2^{-\nu(1+n/2)}\, |s|^{-1}\Big( 2^\nu\, \| f\|_{L_2(K_\nu)} + \| \nabla u\|_{L_2(K_\nu)}+ \| p\|_{L_2(K_\nu)}\Big).
\]
Thus, the estimate (\ref{3l3}) implies (\ref{1l3}). \hfill $\Box$ \\

Using the estimate (\ref{1l3}), one can easily prove the following lemma.

\begin{Le} \label{l4}
Suppose that $u\in W_2^2({\cal D})$, $p\in W_2^1({\cal D})$ for every bounded subdomain ${\cal D}\subset K$, $\overline{\cal D}
\subset \overline{K}\backslash \{ 0\}$, and that $u\in V_{\beta-1}^1(K)$, $p\in V_{\beta-1}^0(K)$. If $(u,p)$
satisfies {\em (\ref{par3})} with $|s|=1$ and with data $f\in V_\beta^0(K)$, $g\in X_\beta^1(K)$, then
$u\in E_\beta^2(K)$, $p\in V_\beta^1(K)$ and
\[
\| u\|_{E_\beta^2(K)} + \| p\|_{V_\beta^1(K)} \le c\, \Big( \| f\|_{V_\beta^0(K)} + \| g\|_{X_\beta^1(K)}
   + \| u\|_{V_{\beta-1}^1(K)} + \| p\|_{V_{\beta-1}^0(K)}\Big).
\]
\end{Le}

P r o o f.
Multiplying (\ref{1l3}) by $2^{2\nu\beta}$, we get
\begin{equation} \label{1l4}
\| \zeta_\nu u\|^2_{E_\beta^2(K)} + \| \zeta_\nu p\|^2_{V_\beta^1(K)} \le
  c\, \Big( \| r^\beta f\|^2_{L_2(K_\nu)} + \| \zeta_\nu g\|^2_{X_\beta^1(K)}
 + \  \| \eta_\nu u \|^2_{V_{\beta-1}^1(K)} +  \| \eta_\nu p\|^2_{V_{\beta-1}^0(K_\nu)}\Big)
\end{equation}
for every $\nu$, where $\eta_\nu=\zeta_{\nu-1}+\zeta_\nu+\zeta_{\nu+1}=1$ on $K_\nu$. The constant $c$ in (\ref{1l4})
is independent of $\nu$. Using the inequalities
\begin{eqnarray*}
&& \sum_{\nu=-\infty}^{+\infty} \Big( \| r^\beta f\|^2_{L_2(K_\nu)} + \| \zeta_\nu g\|^2_{X_\beta^1(K)}
   + \| \eta_\nu u \|^2_{V_{\beta-1}^1(K)} +  \| \eta_\nu p\|^2_{V_\beta^0(K_\nu)}\Big) \\
&& \le c\, \Big( \| f\|^2_{V_\beta^0(K)} + \| g\|^2_{X_\beta^1(K)} + \| u\|^2_{V_{\beta-1}^1(K)} + \| p\|^2_{V_{\beta-1}^0(K)}\Big)
\end{eqnarray*}
and
\begin{eqnarray*}
\| u\|^2_{E_\beta^2(K)} + \| p\|^2_{V_\beta^1(K)} \le c\, \sum_{\nu=-\infty}^{+\infty} \Big( \| \zeta_\nu u\|^2_{E_\beta^2(K)}
  + \| \zeta_\nu p\|^2_{V_\beta^1(K)}\Big),
\end{eqnarray*}
we obtain the desired result. \hfill $\Box$  \\

\subsection{An auxiliary problem for the function \boldmath $p$}

Suppose that $(u,p) \in E_\beta^2(K)\times V_\beta^1(K)$ is a solution of the problem (\ref{par3}).
Multiplying the equation $su-\Delta u+\nabla p = f$ by $\nabla q$ and integrating, we obtain
\begin{equation} \label{Neumann}
\int_K \nabla p\cdot \nabla q\, dx = \langle F, q\rangle  \quad\mbox{ for all } q\in V_{-\beta}^1(K),
\end{equation}
where
\[
\langle F, q\rangle = \int_K  \big( (f+\Delta u)\cdot \nabla q -sg\, q \big) \, dx.
\]
The following lemma allows us to use another representation of the functional $F$.

\begin{Le} \label{l5a}
Suppose that $u\in E_\beta^2(K)$ and $q\in V_{-\beta}^1(K) \cap V_{1-\gamma}^2(K)$, where
$\beta \le \gamma \le \beta+1/2$. Then
\[
\int_K \Delta u \cdot \nabla q\, dx = - \int_K \nabla g\cdot \nabla q\, dx + \int_{\partial K} \sum_{i,j=1}^n
  \frac{\partial u_i}{\partial x_j}\, \Big( n_j\, \frac{\partial q}{\partial x_i} - n_i\, \frac{\partial q}{\partial x_j}
  \Big)\, d\sigma,
\]
where $g=-\nabla\cdot u$.
\end{Le}

P r o o f.
Suppose that $q\in C_0^\infty(\overline{K}\backslash\{ 0\})$. Then
\begin{eqnarray*}
&& \hspace{-1em} \int_K \big( \Delta u + \nabla g\big)\cdot \nabla q\, dx =  \int_K \Big( \sum_{i=1}^3 \Delta u_i \,
  \frac{\partial q}{\partial x_i}-\sum_{j=1}^3 \frac{\partial(\nabla\cdot u)}{\partial x_j}\, \frac{\partial q}{\partial x_j}\Big)\, dx \\
&& = \int_K \sum_{i,j=1}^3 \Big( \frac{\partial^2 u_i}{\partial x_j^2} \, \frac{\partial q}{\partial x_i}
  - \frac{\partial^2 u_i}{\partial x_i\, \partial x_j}\, \frac{\partial q}{\partial x_j}\Big)\, dx
  = \int_K \sum_{i,j=1}^3 \Big( \frac{\partial}{\partial x_j}\, \Big( \frac{\partial u_i}{\partial x_j} \, \frac{\partial q}{\partial x_i}
  \Big) - \frac{\partial}{\partial x_i}\, \Big( \frac{\partial u_i}{\partial x_j}\, \frac{\partial q}{\partial x_j}\Big)\Big)\, dx\\
&& = \int_{\partial K} \sum_{i,j=1}^3  \frac{\partial u_i}{\partial x_j}\, \Big( n_j\, \frac{\partial q}{\partial x_i}
  - n_i\, \frac{\partial q}{\partial x_j} \Big)\, d\sigma.
\end{eqnarray*}
Using the inequalities
\[
\int_{\partial K} r^{-2\gamma+1} \, \big| \nabla q\big|^2\, d\sigma \le  \| \nabla q\|^2_{V_{1-\gamma}^{1/2}(\partial K)}
  \le c\, \| q\|^2_{V_{1-\gamma}^2(K)}
\]
and
\[
\int_{\partial K} r^{2\gamma-1} \, \Big| \frac{\partial u}{\partial x_j}\Big|^2\, d\sigma
  \le \int_{\partial K} \big( r^{2\beta-1}+ r^{2\beta}\big) \, \Big| \frac{\partial u}{\partial x_j}\Big|^2\, d\sigma
  \le c\, \Big\|  \frac{\partial u}{\partial x_j}\Big\|^2_{E_\beta^1(K)}
\]
(cf. (\ref{estV}) and (\ref{estE})), we can easily show that the above identity holds also
for $q\in V_{-\beta}^1(K) \cap V_{1-\gamma}^2(K)$. \hfill $\Box$ \\

We employ a regularity result for solutions of the problem (\ref{Neumann}). For this, we have to consider
the following operator pencil ${\cal N}(\lambda)$ which is defined as
\[
{\cal N}(\lambda)\, U = \Big( -\delta U-\lambda(\lambda+n-2)U\, , \, \frac{\partial U}{\partial\vec{n}}\big|_{\partial\Omega}\Big)\
  \quad\mbox{for } U\in W^2(\Omega).
\]
As is known (see e.g. \cite[Section 2.3]{kmr-01}), the eigenvalues of this pencil are real and
generalized eigenfunctions do not exist. The spectrum contains, in particular, the simple eigenvalues $\mu_1^+=0$ and $\mu_1^-=-1$ with the
eigenfunction $\phi_1 = const.$ The interval $(-1,0)$ is free of eigenvalues. Let $\mu_2^+\le \mu_3^+\le \cdots$ be the
positive and $\mu_j^-=-1-\mu_j^+$ the negative eigenvalues of the pencil ${\cal N}$.

\begin{Le} \label{l5}
{\em 1)} Suppose that $F\in \big( V_{-\beta}^1(K)\big)^*$ and that $-\beta -1/2$ is not an eigenvalues of the pencil
${\cal N}(\lambda)$. Then there exists a unique solution $p\in V_\beta^1(K)$ of the problem {\em (\ref{Neumann})}.

{\em 2)}
Suppose that $p \in V_\beta^1(K)$ is a solution of the problem {\em (\ref{Neumann})}, where $F\in \big( V_{-\beta}^1(K)\big)^*
\cap \big( V_{1-\gamma}^2(K)\big)^*$, $\beta<\gamma$. We assume that the numbers $-\beta-1/2$ and $-\gamma-1/2$
are not eigenvalues of the pencil ${\cal N}(\lambda)$. Then
\[
p = \sum c_j^\pm \, r^{\mu_j^\pm} \phi_j(\omega) + v
\]
where $\mu_j^\pm$ are the eigenvalues of the pencil ${\cal N}(\lambda)$ in the interval $-\gamma-1/2 < \lambda < -\beta-1/2$,
$\phi_j$ are corresponding eigenfunctions, $v\in V_{\gamma-1}^0(K)$ and
\[
\| v\|_{V_{\gamma-1}^0(K)} \le c\, \| F\|_{\big( V_{1-\gamma}^2(K)\big)^*}
\]
with a constant $c$ independent of $F$.
\end{Le}

P r o o f.
The first assertion can be found e.g. in \cite[Theorem 7.3.7]{mr-10}. We prove the second assertion.
Under the conditions of the lemma, the operator
\[
V_{1-\gamma}^2(K) \ni u \to \Big( -\Delta u, \, \frac{\partial u}{\partial n}\Big)
  \in V_{1-\gamma}^0(K) \times V_{1-\gamma}^{1/2}(\partial K)
\]
is an isomorphism (see, e.g., \cite[Theorem 6.1.1]{kmr-97}). Consequently, the adjoint operator is also an isomorphism.
This means that for arbitrary $F\in \big( V_{1-\gamma}^2(K)\big)^*$ there exists a pair
$(v,\phi)\in V_{\gamma-1}^0(K)\times V_{\gamma-1}^{-1/2}(\partial K)$ such that
\[
- \int_K v\, \Delta q\, dx + \int_{\partial K} \phi \, \frac{\partial q}{\partial n}\, d\sigma = \langle F,q\rangle
   \ \mbox{ for all }q\in V_{1-\gamma}^2(K).
\]
This together with (\ref{Neumann}) yields
\[
- \int_K (v-p)\, \Delta q\, dx + \int_{\partial K} (\phi-p)\, \frac{\partial q}{\partial n}\, d\sigma =0
  \ \mbox{ for all }q\in V_{1-\beta}^2(K)\cap V_{1-\gamma}^2(K).
\]
Thus, the triple $(v-p,\phi-p|_{\partial K},0)$ can be understood as a generalized solution of
the Neumann problem
\[
-\Delta(v-p)=0 \ \mbox{ in K}, \quad \frac{\partial(v-p)}{\partial n} =0 \ \mbox{ on }\partial K \backslash \{ 0\}
\]
in the sense of \cite[Section 3.2]{kmr-97}. From \cite[Lemma 3.2.4]{kmr-97} we conclude that
$\chi (v-p)\in W^2(K)$ for every $\chi \in C_0^\infty(\overline{K}\backslash \{ 0\})$. Let $\zeta=\zeta(r)$ be
an infinitely differentiable function on $(0,\infty)$ such that $\zeta(r)=1$ for $0<r<1$, $\zeta(r)=0$ for $r>2$,
and let $\eta=1-\zeta$. With the notation $r=|x|$, the functions $\zeta$ and $\eta$ can be also considered as functions on $K$.
Obviously $\zeta(v-p) \in V_{\beta-1}^0(K)$ and
$\eta(v-p)\in V_{\gamma-1}^0(K)$. Applying \cite[Lemma 6.3.1]{kmr-97}, we obtain
$\zeta(v-p)\in V_{\beta+1}^2(K)$ and $\eta(v-p)\in V_{\gamma+1}^2(K)$. Furthermore,
\[
\Delta(\zeta p - \zeta v) = - \Delta(\eta p-\eta v) \in V_{\beta+1}^0(K)\cap V_{\gamma+1}^0(K), \quad
  \frac{\partial(\zeta p -\zeta v)}{\partial n} = \frac{\partial(\eta p -\eta v)}{\partial n}=0 \ \mbox{ on }\partial K
  \backslash\{ 0\}.
\]
Using \cite[Theorem 6.1.4]{kmr-97}, we get
\[
\zeta(p-v) = \sum c_j^\pm \, r^{\mu_j^\pm} \phi_j(\omega) + \eta(v-p)
\]
This proves the lemma. \hfill $\Box$

\subsection{Normal solvability of the mapping \boldmath $A_\beta$}

\begin{Th} \label{t2}
Suppose that $\mbox{\em Re}\, s \ge 0$, $|s|=1$, that the line $\mbox{\em Re}\, \lambda = -\beta+1/2$ does not contain
eigenvalues of the pencil ${\cal L}(\lambda)$, and $-\beta-1/2$ is not an eigenvalue of the pencil ${\cal N}(\lambda)$.
Then the range of the operator {\em (\ref{1t2})} is closed and the kernel has finite dimension.
\end{Th}

P r o o f.
It suffices to show that every $(u,p)\in \big( E_\beta^2(K) \cap \stackrel{\circ}{E}\!{}_{\beta}^1(K)\big)
\times V_\beta^1(K)$ satisfies the estimate
\[
\| u\|^2_{E_\beta^2(K)} + \| p\|^2_{V_\beta^1(K)} \le c\, \Big( \| f\|^2_{E_\beta^0(K)}
  + \| g \|^2_{X_\beta^1(K)} + \int\limits_{\substack{K \\ c_1 < |x| < c_2}}
  \big( |u|^2 + |\nabla u|^2 + |p|^2\big)\, dx\Big)
\]
with certain positive $c_1$ and $c_2>c_1$.
Multiplying (\ref{1l3}) by $2^{2\nu\beta}$ and summing up over all integer $\nu\ge N$,  we get
\begin{eqnarray} \label{2t2}
&&\int\limits_{\substack{K \\ |x|>2^N}} r^{2\beta}\, \Big( |D^2 u|^2 +  |u|^2 + |\nabla p|^2\Big)\, dx
  \le c\, \int\limits_{\substack{K \\ |x|>2^{N-1}}} r^{2\beta}\, \Big( |f|^2 + |\nabla g|^2 + r^{-2}\, |g|^2 \Big)\, dx\nonumber \\
&& \qquad + \ c\, \Big(\int\limits_{\substack{K \\ |x|>2^{N-1}}} r^{2\beta-2}\, \big(  |\nabla u|^2 + r^{-2}\, |u|^2 + |p|^2\big)\, dx
 + \| \eta_N g\|^2_{(V_{-\beta}^1(K))^*}\Big) ,
\end{eqnarray}
where $\eta_N = \zeta_N + \zeta_{N+1}+\cdots$. Here, the constant $c$ is independent of $u,p$ and $s$.

Let $\zeta$ be a two times continuously differentiable function on $\overline{K}$,
$\zeta(x)=1$ for $|x|<\epsilon$  and $\zeta(x)=0$ for $|x|>2\epsilon$. Then
\[
-\Delta(\zeta u_i) + \partial_{x_i} (\zeta p) = \zeta f_i -s\zeta u_i-2\nabla\zeta\cdot \nabla u_i-u_i\, \Delta\zeta
  + p\, \partial_{x_i}\zeta =: F_i
\]
for $i=1,2,3$ and
\[
-\nabla\cdot(\zeta u) = \zeta g - u\cdot \nabla \zeta.
\]
Since the line $\mbox{Re}\, \lambda = 2-\beta-n/2$ is free of eigenvalues of the pencil ${\cal L}(\lambda)$,
there exists a constant $c$ such that
\[
\| \zeta u\|_{V_\beta^2(K)} + \| \zeta p\|_{V_\beta^1(K)} \le c\, \Big( \| F\|_{V_\beta^0(K)} +
  \| \zeta g - u\cdot\nabla\zeta\|_{V_\beta^1(K)}\Big),
\]
where $c$ is independent of $\zeta$. If $\epsilon$ is sufficiently small, then it follows that
\[
\| \zeta u\|^2_{E_\beta^2(K)} + \| \zeta p\|^2_{V_\beta^1(K)} \le c\, \Big( \| \zeta f\|^2_{V_\beta^0(K)}+\| \zeta g\|^2_{V_\beta^1(K)}\Big)
  + C(\epsilon) \, \int\limits_{\substack{K \\ \epsilon < |x|<2\epsilon}} \big(|u|^2+|\nabla u|^2+ |p|^2\big)\, dx.
\]
Combining this with (\ref{2t2}), we obtain
\begin{eqnarray} \label{4t2}
\|  u\|^2_{E_\beta^2(K)} + \| p\|^2_{V_\beta^1(K)} & \le & c\, \Big( \| f\|^2_{V_\beta^0(K)} + \| g\|^2_{V_\beta^1(K)}
  + \| g\|^2_{(V_{-\beta}^1(K))^*} \nonumber \\
&& + \ \int\limits_{\substack{K \\ |x|>\epsilon}} r^{2\beta-2}\, \Big(   |\nabla u|^2 + r^{-2}\, |u|^2  + |p|^2\Big)\, dx \Big).
\end{eqnarray}
Obviously,
\begin{equation} \label{5t2}
\int\limits_{\substack{K \\ |x|>N}} r^{2\beta}\, \Big( r^{-2}\, |\nabla u|^2 + r^{-4}\, |u|^2\Big)\, dx \le \frac{c}{N^2} \,
  \| u\|^2_{E_\beta^2(K)}
\end{equation}
We estimate the integral of $r^{2\beta-2}|p|^2$. By (\ref{Neumann}), there is the decomposition $p= p_1+p_2$, where $p_j$ are the solutions
of the problems
\[
\int_K \nabla p_j \cdot \nabla q\, dx = \langle \Phi_j,q \rangle \ \mbox{ for all }q\in V_{-\beta}^1(K)
\]
with the functionals
\[
\langle \Phi_1,q\rangle = \int_K \big( (f-\nabla g)\cdot\nabla q -sg\, q\big)\, dx, \quad
\langle \Phi_2,q\rangle = \int_K \big( \Delta u + \nabla g\big)\cdot\nabla q \, dx.
\]
Obviously, $\Phi_1 \in \big( V_{-\beta}^1(K)\big)^*$ and $\Phi_2 \in \big(V_{-\beta}^1(K)\big)^*$,
\[
\| \Phi_1\|_{(V_{-\beta}^1(K))^*} \le  \| f\|_{V_\beta^0(K)} + \| \nabla g\|_{V_\beta^0(K)} + \| g\|_{(V_{-\beta}^1(K))^*}\, .
\]
Moreover, by Lemma \ref{l5a}, $\Phi_2\in (V_{1-\gamma}^2(K))^*$ for $\beta \le \gamma \le \beta+1/2$ and
\begin{eqnarray*}
\big|\langle \Phi_2,q\rangle\big|^2 & = & \Big| \int_{\partial K} \sum_{i,j=1}^n
  \frac{\partial u_i}{\partial x_j}\, \Big( n_j\, \frac{\partial q}{\partial x_i} - n_i\, \frac{\partial q}{\partial x_j}
  \Big)\, d\sigma\Big|^2 \\
& \le & c\, \int_{\partial K} r^{2\gamma-1} \, \Big| \frac{\partial u}{\partial x_j}\Big|^2\, d\sigma\
  \int_{\partial K} r^{1-2\gamma} \, \big| \nabla q\big|^2\, d\sigma \le c\, \| u\|^2_{E_\beta^2(K)}\ \| q\|^2_{V_{1-\gamma}^2(K)}\, .
\end{eqnarray*}
For the last inequality, we used (\ref{estV}) and (\ref{estE}). By Lemma \ref{l5},
\[
\| p_1 \|^2_{V_\beta^1(K)} \le c_1\, \Big( \| f\|^2_{V_\beta^0(K)} + \| \nabla g\|^2_{V_\beta^0(K)}
  + \| g\|^2_{(V_{-\beta}^1(K))^*}\Big).
\]
Let $\gamma \in (\beta,\beta+1/2]$ be such that the interval $-\gamma-1/2 \le \lambda \le -\beta-1/2$
is free of eigenvalues of the pencil ${\cal N}(\lambda)$. Then, by the second part of Lemma \ref{l5},
$p_2 \in V_{\gamma-1}^0(K)$ and
\[
\| p_2\|_{V_{\gamma-1}^0(K)} \le c\, \| \Phi_2 \|_{(V_{1-\gamma}^2(K))^*} \le c_2\, \| u\|_{E_\beta^2(K)}.
\]
Consequently,
\begin{eqnarray*}
&& \int\limits_{\substack{K \\ |x|>N}} r^{2\beta-2}\, |p|^2\, dx \le 2\, \| p_1\|^2_{V_{\beta-1}^0(K)}
  + \frac{2}{N^{2(\gamma-\beta)}} \, \| p_2\|^2_{V_{\gamma-1}^0(K)} \\
&&  \le 2c_1\,\Big( \| f\|^2_{V_\beta^0(K)} + \| \nabla g\|^2_{V_\beta^0(K)}
  + \| g\|^2_{(V_{-\beta}^1(K))^*}\Big)+ \frac{c_2^2}{N^{2(\gamma-\beta)}} \, \| u\|^2_{E_\beta^2(K)}.
\end{eqnarray*}
This together with (\ref{4t2}) implies
\begin{eqnarray*}
\|  u\|^2_{E_\beta^2(K)} + \| p\|^2_{V_\beta^1(K)} & \le & c\, \Big( \| f\|^2_{V_\beta^0(K)} + \| g\|^2_{V_\beta^1(K)}
  +  \| g\|^2_{(V_{-\beta}^1(K))^*}\Big) \\
&& +\ C(\epsilon,N) \int\limits_{\substack{K \\ \epsilon<|x|<N}} \Big( |\nabla u|^2 + |u|^2
  + |p|^2\Big)\, dx
\end{eqnarray*}
if $N$ is sufficiently large. The proof is complete. \hfill $\Box$ \\

The following two lemmas show that the conditions of the last theorem are necessary.

\begin{Le} \label{l6b}
The assertion of Theorem {\em \ref{t2}} is not true if the line $\mbox{\em Re}\, \lambda = -\beta+1/2$
contains an eigenvalue of the pencil ${\cal L}(\lambda)$.
\end{Le}

P r o o f.
Let $\lambda$ be an eigenvalue of the pencil ${\cal L}$, $\mbox{Re}\, \lambda = -\beta+1/2$, and let $(u_0,p_0)$
be an eigenvector corresponding to this eigenvalue. Then the pair $(v,q)=\big(r^\lambda u_0(\omega),\, r^{\lambda-1}p_0(\omega)\big)$ is a solution of
the system
\[
-\Delta v+\nabla q =0,  \ \  \nabla\cdot v =0 \ \mbox{ in }K
\]
satisfying the condition $v=0$ on $\partial K \backslash \{ 0\}$. We define by $\zeta_\varepsilon$ a smooth function
on ${\Bbb R}_+$ such that $\zeta_\varepsilon(r)=1$ for $\varepsilon<r<1$, $\zeta_\varepsilon=0$ outside the interval
$(\varepsilon/2,2)$ and $|\zeta_\varepsilon^{(j)}(r)| \le c\, r^{-j}$ for $j\le 2$, where the constant $c$ is
independent of $r$ and $\varepsilon$. Obviously, the functions $u(x)=\zeta_\varepsilon(r)\, v(x)$ and
$p(x)=\zeta_\varepsilon(r)\, q(x)$ satisfy the estimate
\begin{eqnarray*}
\| u\|^2_{E_\beta^2(K)} + \| p\|^2_{V_\beta^1(K)} & \ge & \int_K \Big( r^{2\beta-4}\,  |u|^2+ r^{2\beta-2}\, |p|^2\Big)\, dx
  \ge c\int_\varepsilon^1 r^{2\beta+2\mbox{\scriptsize Re}\lambda-2}\, dr = -c\log\varepsilon.
\end{eqnarray*}
Obviously, $|u|\le c\, r^{\mbox{\scriptsize Re}\lambda}$ and
$|\nabla p - \Delta u|\le c\, r^{\mbox{\scriptsize Re}\lambda-2}$. Furthermore, $\nabla p -\Delta u$
vanishes outside the sets $\{ x:\, \varepsilon/2 < |x|<\varepsilon\}$ and $\{ x:\, 1<|x|<2\}$. Thus,
\begin{eqnarray*}
\| su-\Delta u + \nabla p \|^2_{E_\beta^0(K)} &\le & c\, \Big( \int_0^2 r^{2\beta+2\mbox{\scriptsize Re}\lambda+2}\, dr
  + \int_{\varepsilon/2}^\varepsilon   r^{2\beta+2\mbox{\scriptsize Re}\lambda-2}\, dr
  + \int_1^2 r^{2\beta+2\mbox{\scriptsize Re}\lambda-2}\, dr\Big) \\
& = & 2c(2+\log 2).
\end{eqnarray*}
Using the estimates $|\nabla\cdot u|\le c\, r^{\mbox{\scriptsize Re}\lambda-1}$,
$\big|\nabla(\nabla\cdot u)\big|\le c\, r^{\mbox{\scriptsize Re}\lambda-2}$ and the fact that
$\nabla\cdot u=v\cdot\nabla\zeta_\varepsilon$ vanishes outside the sets $\{ x:\, \varepsilon/2 < |x|<\varepsilon\}$
and $\{ x:\, 1<|x|<2\}$, we analogously obtain
\[
\| \nabla \cdot u\|_{V_\beta^1(K)} \le C \ \mbox{ and } \ \| \nabla \cdot u\|_{(V_{-\beta}^1(K))^*}
  \le \| \nabla \cdot u\|_{V_{\beta+1}^0(K)} \le C.
\]
The space $E_\beta^2(K)$ is compactly imbedded into $V_{\beta-1}^0(K)$, while $V_\beta^1(K)$ is compactly imbedded
into the space $V_{\beta,\beta-2}^0(K)$ with the norm
\[
\| p\|_{V_{\beta,\beta-2}^0(K)} = \Big( \int\limits_{\substack{K \\ |x|<1}} r^{2\beta}\, |p|^2\, dx +
  \int\limits_{\substack{K \\ |x|>1}} r^{2\beta-4}\, |p|^2\, dx\Big)^{1/2}.
\]
One can easily check that the $V_{\beta-1}^0(K)$-norm of $u$ and the ${V_{\beta,\beta-2}^0(K)}$ of $p$ have
an upper bound $C$ which is independent of $\varepsilon$. Thus, there does not exist a constant $c$ independent
of $\varepsilon$ such that
\begin{eqnarray} \label{1l6b}
\| u\|_{E_\beta^2(K)} + \| p\|_{V_\beta^1(K)} & \le & c\, \Big( \| su-\Delta u + \nabla p\|_{V_\beta^0(K)}
  + \| \nabla\cdot u\|_{X_\beta^1(K)} \nonumber \\ && \quad
  + \ \| u\|_{V_{\beta-1}^0(K)}  + \| p\|_{V_{\beta,\beta-2}^0(K)}\Big).
\end{eqnarray}
The result follows. \hfill $\Box$ \\

We prove the same result concerning the condition on the eigenvalues of the pencil ${\cal N}(\lambda)$.

\begin{Le} \label{l6a}
The assertion of Theorem {\em \ref{t2}} is not true if $-\beta-1/2$ is an eigenvalue of the pencil ${\cal N}(\lambda)$.
\end{Le}

P r o o f.
It suffices to show that there does not exist a constant $c$ such that the estimate (\ref{1l6b}) is valid for all
$u\in  E_\beta^2(K)\cap \stackrel{\circ}{E}\!{}_\beta^1(K)$ and $p\in V_\beta^1(K)$ if $-\beta-1/2$
is an eigenvalue of the pencil ${\cal N}(\lambda)$.

Let $\phi$ be an eigenfunction corresponding to the eigenvalue $\lambda=-\beta-1/2$ of the pencil ${\cal N}$. Then the function
$p_0(x)=r^\lambda \, \phi(\omega)$ is a solution of the problem
\[
\Delta p_0 =0\ \mbox{ in }K,\quad \frac{\partial p_0}{\partial n}=0\ \mbox{ on }\partial K\backslash\{ 0\}.
\]
Obviously, $\nabla p_0=r^{\lambda-1}\, \psi(\omega)$ with a certain smooth function $\psi$ on $\Omega$.
We define $v_0= -s^{-1}\, \nabla p_0$. Then
\[
(s-\Delta) v_0+\nabla p_0=0, \quad  \nabla\cdot v_0 =0 \ \mbox{ in }K
\]
and $v_0 \cdot n=0$ on $\partial K\backslash \{ 0\}$. We correct the function $v_0$ by a term with support
in a neighborhood of the boundary. Let $\nu(x)$ denote the distance of $x$ from
the boundary $\partial K$. Obviously the function $\nu$ is positively homogeneous of degree 1. Furthermore, it is smooth
in the neighborhood $\nu <\delta|x|$ of the boundary $\partial K$, where $\delta$ is a sufficiently small positive number.
Moreover, $|\nabla\nu|=1$ for $\nu<\delta|x|$ and $\nabla\nu(x)$ is orthogonal to $\partial K$ for $x\in \partial K$.

We introduce the functions
\[
v_{0,\nu} = v_0\cdot \nabla\nu \ \mbox{ and } \ v_{0,\tau} = v_0-v_{0,\nu}\nabla\nu
\]
in the neighborhood $\nu<\delta|x|$ of $\partial K$ and define
\[
u_0(x) = v_0(x) - \chi\Big( \frac{\nu}r \Big) \, e^{-\nu\sqrt{s}}\, v_{0,\tau}(x),
\]
where $\chi$ is a smooth cut-off function, $\chi=1$ in $(0,\delta/2)$, $\chi=0$ in $(\delta,\infty)$.
For $\nu>\delta|x|$, we set $u_0(x)=v_0(x)$. Here by $\sqrt{s}$, we mean the square root of $s$
with positive real part. Then $u_0=v_0-v_{0,\tau}=v_{0,\nu}\nabla\nu =0$ on $\partial K$.
Since $(s-\Delta)\, e^{-\nu\sqrt{s}} = \sqrt{s}\, e^{-\nu\sqrt{s}} \, \Delta\nu$, we obtain
\begin{equation} \label{2l6a}
(s-\Delta) u_0 + \nabla p_0 = - (s-\Delta)\, \chi\Big( \frac{\nu}r \Big) \, e^{-\nu\sqrt{s}}\, v_{0,\tau}
  =  \chi\Big( \frac{\nu}r \Big) \, e^{-\nu\sqrt{s}} \, f_0   + F_0,
\end{equation}
where
\[
f_0=-\sqrt{s}\, v_{0,\tau}\Delta\nu- 2\sqrt{s}\sum_{j=1}^3 \frac{\partial\nu}{\partial x_j}\, \frac{\partial v_{0,\tau}}{\partial x_j}
  + \Delta v_{0,\tau} , \quad F_0 = \Big[ \Delta,\chi\Big(\frac{\nu}r\Big)\Big]\, e^{-\nu\sqrt{s}}\, v_{0,\tau}
\]
(here $[\Delta,\chi]=\Delta\chi-\chi\Delta$ denotes the commutator of $\Delta$ and $\chi$). Furthermore,
\begin{equation} \label{1l6a}
\nabla\cdot u_0 = - \nabla\cdot \chi\Big( \frac{\nu}r \Big) \, e^{-\nu\sqrt{s}}\, v_{0,\tau}
  = - \chi\Big( \frac{\nu}r \Big) \, e^{-\nu\sqrt{s}}\, \nabla \cdot v_{0,\tau}
  - e^{-\nu\sqrt{s}}\, v_{0,\tau}\cdot \nabla \chi\Big( \frac{\nu}r \Big)
\end{equation}
since $v_{0,\tau}\cdot \nabla e^{-\nu\sqrt{s}}= -\sqrt{s}\, e^{-\nu\sqrt{s}}\, v_{0,\tau}\cdot\nabla\nu=0$.
Let $\zeta_N$ be a two times continuously differentiable function on ${\Bbb R}_+=(0,\infty)$ such that
\[
\zeta_N(r)=0 \ \mbox{ for }r<1 \mbox{ and }r>2N, \quad \zeta_N(r)=1 \ \mbox{ for }2<r<N.
\]
We may assume that $|r^k\zeta_N^{(k)}(r)|\le c$ for $k=0,1,2$ with a constant $c$ independent of $N$.
We consider the functions
\[
u = \zeta_N u_0 \ \mbox{ and } p= \zeta_N p_0
\]
Obviously, $u \in E_\beta^2(K)$, $p\in V_\beta^1(K)$ and $u=0$ on $\partial K$. It can be easily checked that
\[
\| u\|_{E_\beta^2(K)}+\| p\|_{V_\beta^1(K)}
  \ge \| \zeta_N(r)\, p_0\|_{V_{\beta-1}^0(K)} \ge C\, \log N
\]
with a positive constant $C$. We estimate the $E_\beta^0(K)$-norm of $(s-\Delta)u + \nabla p$.
Obviously, $(s-\Delta)u+\nabla p =F_1+F_2$, where
\[
F_1 = \zeta_N\big( (s-\Delta)u_0+\nabla p_0\big), \ \mbox{ and }\
F_2=p_0\nabla\zeta_N-2\sum_{j=1}^3 \frac{\partial\zeta_N}{\partial x_j}\, \frac{\partial u_0}{\partial x_j} -u_0\Delta\zeta_N.
\]
Here, $|F_1|\le c\, r^{\lambda-2}$ and $|F_2|\le c\, r^{\lambda-1}$. Since $F_2$ is zero outside the regions $1<|x|<2$
and $N<|x|<2N$, we get
\[
\| (s-\Delta)u+\nabla p\|^2_{E_\beta^0(K)} \le c\, \Big( \int_1^{2N} r^{-3}\, dr + \int_1^2 r^{-1}\, dr + \int_N^{2N} r^{-1}\, dr\Big)
  < \frac c2 + 2c\log 2.
\]
Next, we estimate the $X_\beta^1(K)$-norm of $\nabla\cdot u$. Since $|\partial_x^\alpha \nabla\cdot u_0|\le c\, r^{\lambda-2}$
for $|\alpha|\le 1$ and $r>1$, we have $|\partial_x^\alpha \nabla\cdot u|\le c\, r^{\lambda-2}$ for $|\alpha|\le 1$. Thus,
\[
\| \nabla\cdot u\|^2_{V_\beta^1(K)} \le c\, \int_1^{2N} r^{-3}\, dr < \frac c2\, .
\]
Since $|u_0|\le c\, r^{\lambda-1}$ und $|\nabla \zeta_N|\le c\, r^{-1}$, we get
\[
\| u_0\cdot \nabla\zeta_N \|^2_{(V_{-\beta}^1(K))^*} \le \| u_0\cdot \nabla\zeta_N \|^2_{V_{\beta+1}^0(K)}
  \le c\, \Big( \int_1^2 r^{-1}\, dr + \int_N^{2N} r^{-1}\, dr\Big) = 2c\log 2.
\]
Furthermore, we conclude from (\ref{1l6a}) that $|\zeta_N\nabla\cdot u_0|\le c\, r^{\lambda-2}\, e^{-\nu\sqrt{s}}$.
Consequently,
\begin{eqnarray*}
\| \zeta_N\nabla\cdot u_0 \|^2_{(V_{-\beta}^1(K))^*} & \le &
\|  \zeta_N\nabla\cdot u_0 \|^2_{V_{\beta+1}^0(K)}
\le c\, \int r^{-3}\, e^{-2\nu(x)\, \mbox{\scriptsize Re}\, \sqrt{s}}\, dx,
\end{eqnarray*}
where the integration is extended over the set of all $x\in K$ such that $1<r=|x|<2N$ and $\nu(x)<\delta r$.
For arbitrary $x\in K$, $\nu(x)<\delta|x|$, let $\tau(x)$ be the nearest point to $x$ on $\partial K$.
Obviously, $|\tau(x)|\le |x|\le (1-\delta)^{-1/2}|\tau(x)|$. This implies
\begin{eqnarray*}
\| \zeta_N\nabla\cdot u_0  \|^2_{(V_{-\beta}^1(K))^*} \le  c\, \int\limits_{\substack{\partial K \\ |\tau|^2>1-\delta}}
  \int\limits_0^{\delta r} |\tau|^{-3}  e^{-2\nu\, \mbox{\scriptsize Re}\, \sqrt{s}}\, d\nu\, d\tau
\le  \frac{c}{2\, \mbox{\small Re}\, \sqrt{s}} \int\limits_{\substack{\partial K \\ |\tau|^2>1-\delta}} |\tau|^{-3} d\tau
  \le  \frac{c'}{\mbox{\small Re}\, \sqrt{s}}\, .
\end{eqnarray*}
The last estimates yield
\[
\| \nabla\cdot u\|_{X_\beta^1(K)} \le C
\]
with a constant $C$ independent of $N$. Finally, since $|p(x)|\le c\, r^\lambda$, $|u(x)| \le c\, r^{\lambda-1}$
and $u(x),p(x)$ are zero for $|x|<1$, we get the estimate
\[
\| u\|^2_{V_{\beta-1}^0(K)} + \| p\|^2_{V^0_{\beta,\beta-2}(K)} \le c\,  \int\limits_{\substack{K \\ |x|>1}}
  r^{2\beta+2\lambda-4}\, dx = C < \infty.
\]
Thus, the right-hand side of (\ref{1l6b}) has an upper bound independent of $N$ for our pair $(u,p)$, while
the left-hand side can be estimate from below by $C\, \log N$ with positive $C$. This proves the lemma. \hfill $\Box$

\subsection{A regularity result for solutions of the problem (\ref{par3})}

Using Lemma \ref{l5} and regularity results for solutions of the stationary Stokes system,
we can prove the following lemma.

\begin{Le}  \label{l5b}
Suppose that $(u,p) \in E_\beta^2(K) \times V_\beta^1(K)$ is a solution of the problem {\em (\ref{par3})}, where
\[
f\in E_\beta^0(K) \cap E_\gamma^0(K), \quad g\in X_\beta^1(K)\cap X_\gamma^1(K).
\]
We assume that one of the following two conditions is satisfied:
\begin{itemize}
\item[{\em (i)}] $\beta<\gamma$ and the interval $-\gamma-1/2 \le \lambda \le -\beta-1/2$ does not contain eigenvalues of the pencil
${\cal N}(\lambda)$,
\item[{\em (ii)}] $\beta>\gamma$ and the strip $-\beta+1/2 \le \mbox{\em Re}\, \lambda\le -\gamma+1/2$ is free of eigenvalues
  of the pencil ${\cal L}(\lambda)$.
\end{itemize}
Then $u\in E_\gamma^2(K)$, $p\in V_\gamma^1(K)$ and
\[
\| u\|_{E_\gamma^2(K)} +\| p\|_{V_\gamma^1(K)} \le c\, \Big( \| f\|_{E_\gamma^0(K)} + \| g\|_{X_\gamma^1(K)}
  + \| u\|_{E_\beta^2(K)} +\| p\|_{V_\beta^1(K)}\Big).
\]
\end{Le}

P r o o f.
1) Suppose that the condition (i) is satisfied. In addition, let $\gamma\le \beta+1/2$. Then $p$ is a
solution of the problem (\ref{Neumann}), where $F\in \big(V_{-\beta}^1(K)\big)^* \cap \big( V_{1-\gamma}^2(K)\big)^*$,
\begin{eqnarray} \label{1l5b}
\langle F,q\rangle & = & \int_K \big( (f+\Delta u)\cdot \nabla q -sg\, q\big)\, dx \ \mbox{ for }q\in V_{-\beta}^1(K) \quad\mbox{and}\\
\langle F,q\rangle & = &  \int_K \Big( (f-\nabla g)\cdot \nabla q-sg\, q\Big)\, dx + \int_{\partial K} \sum_{i,j=1}^n
  \frac{\partial u_i}{\partial x_j}\, \Big( n_j\, \frac{\partial q}{\partial x_i} - n_i\, \frac{\partial q}{\partial x_j}
  \Big)\, d\sigma \label{2l5b}
\end{eqnarray}
for $q\in V_{1-\gamma}^2(K)$ (cf. Lemma \ref{l5a}). From Lemma \ref{l5} we conclude that
$p\in V_{\gamma-1}^0(K)$.
Furthermore, $u\in E_\beta^2(K) \subset V_{\beta-1}^1(K) \cap V_\beta^1(K) \subset V_{\gamma-1}^1(K)$.
Applying Lemma \ref{l4}, we obtain $u\in E_\gamma^2(K)$ and $p\in V_\gamma^1(K)$. This proves the lemma
for the case $\gamma\le \beta+1/2$. By induction, the assertion of the lemma holds for arbitrary
$\gamma >\beta$ provided that the interval $[-\gamma-1/2,-\beta-1/2]$ is free of eigenvalues of the pencil
${\cal N}(\lambda)$.

2) We assume now that the condition (ii) is satisfied.
Let $\zeta$ be a smooth function on $\overline{K}$ with compact support which is equal to one near the vertex of the cone.
Furthermore, let $\eta=1-\zeta$. Obviously, $\eta(u,p)\in E_\gamma^2(K)\times V_\gamma^2(K)$ since $\gamma<\beta$.
Furthermore, the pair $(\zeta u,\zeta p)$ satisfies the equations
\[
-\Delta(\zeta u) + \nabla(\zeta p) = F,\quad -\nabla\cdot (\zeta u)=G \ \mbox{ in }K
\]
and the boundary condition $\zeta u=0$ on $\partial K\backslash\{ 0\}$, where
\[
F=\zeta f -s\zeta u-2\sum_{j=1}^3 \frac{\partial\zeta}{\partial x_j}\, \frac{\partial u}{\partial x_j} -u\, \Delta\zeta+p\, \nabla\zeta,
\quad G=\zeta g -u\cdot\nabla\zeta.
\]
Suppose first that $\beta-2\le \gamma<\beta$. Then $F\in V_\gamma^0(K)$ and $G\in V_\gamma^1(K)$. Since the
strip $-\beta+1/2 \le \mbox{Re}\, \lambda\le -\gamma+1/2$ is free of eigenvalues of the pencil ${\cal L}(\lambda)$
it follows that $\zeta(u,p) \in V_\gamma^2(K)\times V_\gamma^1(K)$ (see, e.~g., \cite[Chapter 3, Theorem 5.6]{np-94}).
By induction, the same result holds for $\gamma<\beta-2$. Thus, $(u,p) \in E_\gamma^2(K)\times V_\gamma^1(K)$. \hfill $\Box$ \\

Note that the result of Lemma \ref{l5b} is also true if $\beta<\gamma$, the interval $-\gamma-1/2
\le \lambda \le -\beta-1/2$ contains only the eigenvalue $\mu_1^- =-1$ as an inner point, and $\int_K g\, dx =0$
(see Lemma \ref{l11} below).

\subsection{Injectivity of the operator \boldmath $A_\beta$}

Let $\lambda_1^+$ be the smallest positive eigenvalue of the pencil ${\cal L}(\lambda)$ and let $\lambda_1^-=-1-\lambda_1^+$.
Then the strip $\lambda_1^- < \mbox{Re}\, \lambda <\lambda_1^+$ is the widest strip containing the line
$\mbox{Re}\,\lambda=-1/2$ which is free of eigenvalues of the pencil ${\cal L}(\lambda)$. As mentioned above,
$\lambda_1^+\le 1$ and $\lambda_1^- \ge -2$.

\begin{Le} \label{l6}
Suppose that $\mbox{\em Re}\, s \ge 0$, $|s|=1$,  $-\mu_2^+ -1/2 < \beta < \lambda_1^+ +3/2$ and $\beta\not=-1/2$.
Then $A_\beta$ is injective.
\end{Le}

P r o o f.
1) If  $0\le \beta\le 1$, then the assertion of the lemma follows from Theorem \ref{t1} and Lemma \ref{l1b}.

2) Suppose that $1<  \beta < \lambda_1^+ +3/2$, i.e., $\lambda_1^- < -\beta+1/2 < -1/2$.
Then the strip  $-\beta+1/2 \le \mbox{Re}\, \lambda\le -1/2$ is free of eigenvalues
of the pencil ${\cal L}(\lambda)$ and it follows from Lemma \ref{l5b} that
$\mbox{ker}\, A_\beta \subset\mbox{ker}\, A_1$. Since $A_1$ is injective, it follows that
the kernel of $A_\beta$ is trivial.

3) If $-1/2 <\beta < 0$, then the interval $-1/2 \le \lambda \le -\beta-1/2$
is free of eigenvalues of the pencil ${\cal N}(\lambda)$ (even the interval $(-1,0)$ is free of eigenvalues
of this pencil). Thus, it follows from Lemma \ref{l5b} that $\mbox{ker}\, A_\beta \subset\mbox{ker}\, A_0$.
This implies that $\mbox{ker}\, A_\beta$ is trivial.

4) We consider the case $\max(-1,-\mu_2^+-1/2) < \beta <-1/2$ and set $\gamma=\beta+1/2$. If $(u,p)\in
\mbox{ker}\, A_\beta$, then
\[
\int_K \nabla p \cdot \nabla q \, dx = \langle F,q\rangle = \int_K \Delta u\cdot \nabla q\, dx \ \mbox{ for all }q\in V_{-\beta}^1(K).
\]
By Lemma \ref{l5a}, the functional $F$ is continuous on $V_{-\beta}^1(K)$ and on $V_{1-\gamma}^2(K)$.
Since $\mu_1^-=-1 < -\gamma-1/2<0<-\beta-1/2 < \mu_2^+$, the interval  $-\gamma-1/2 \le \lambda\le -\beta-1/2$
contains only the eigenvalue $\mu_1^+=0$ of the pencil ${\cal N}(\lambda)$,
and the corresponding eigenfunction is constant. Thus, we get the decomposition
$p=c+q$, where $q\in V_{\gamma-1}^0$ and $c$ is a constant (see Lemma \ref{l5}).
Obviously, $E_\beta^2(K) \subset V_{\gamma-1}^1(K)$. This means that the pair $(u,q)$ belongs to
$V_{\gamma-1}^1(K) \times V_{\gamma-1}^0(K)$. Furthermore,
\[
su-\Delta u + \nabla q =0, \quad \nabla\cdot u=0 \ \mbox{ in }  K, \quad u=0\ \mbox{ on }\partial K\backslash\{ 0\}.
\]
Using Lemma \ref{l4}, we conclude that $(u,q)\in E_\gamma^2(K) \times V_\gamma^1(K)$, i.e., $(u,q) \in \mbox{ker}\,A_\gamma$.
As was shown above, the kernel of $A_\gamma$ is trivial. This means that $u=0$ and $q=0$. The last implies that $p$ is
constant in $K$. However, the space $V_\beta^1(K)$ contains only the constant $p=0$. This shows that $\mbox{ker}\, A_\beta$ is trivial.

5) Suppose that $\mu_2^+ > 1/2$ and $-\mu_2^+-1/2 < \beta \le -1$. Let $\gamma$ be an arbitrary number in the interval
$(-1,-\frac 12 )$. Then $\mu_1^+=0 < -\gamma-1/2 < -\beta-1/2 < \mu_2^+$. Since the interval $-\gamma-1/2 \le \lambda \le -\beta-1/2$
is free of eigenvalues of the pencil ${\cal N}(\lambda)$, it follows from Lemma \ref{l5b} that
$\mbox{ker}\, A_\beta \subset \mbox{ker}\, A_\gamma$. As was shown above, the kernel of $A_\gamma$ and, therefore, also
the kernel of $A_\beta$ are trivial. The proof is complete. \hfill $\Box$

\subsection{Bijectivity of the operator of the problem (\ref{par3})}

We are interested in conditions on $\beta$ which ensure the bijectivity of the operator $A_\beta$.
By Lemma \ref{l6a}, the operator $A_\beta$ is not Fredholm if $-\beta-1/2$ (and, therefore, also $\beta-1/2$)
is an eigenvalue of the pencil ${\cal N}(\lambda)$. In particular, we must exclude the value $\beta=1/2 = \mu_1^+ +1/2$.
We consider the cases $\beta<1/2$ and $\beta>1/2$ separately. \\

\noindent {\bf 1) The case \boldmath $\beta < 1/2$} \\

In order to show that the operator $A_\beta$ is surjective for $-\lambda_1^+ +1/2 < \beta < 1/2$, we prove
the following regularity assertion for weak solutions of the problem (\ref{par3}).

\begin{Le} \label{l7}
Let $|s|=1$, $\mbox{\em Re}\, s \ge 0$, and let $(u,p) \in \stackrel{\circ}{E}\!{}_0^1(K) \times \big( L_2(K) + V_0^1(K)\big)$
be the solution of the problem {\em (\ref{1t1})}, where
\[
f\in E_0^{-1}(K)\cap E_\beta^0(K), \quad g\in L_2(K)\cap \big( V_0^1(K)\big)^* \cap X_{\beta}^1(K)
\]
and $-\lambda_1^+ +1/2 < \beta < 1/2$. Then $(u,p) \in E_\beta^2(K) \times V_\beta^1(K)$.
\end{Le}

P r o o f.
Let $\zeta$ be the same cut-off function as in the proof of Lemma \ref{l5b}, and let $\eta=1-\zeta$.
Obviously, $\zeta(u,p)\in \stackrel{\circ}{V}\!{}_0^1(K) \times L_2(K)$,
\[
-\Delta(\zeta u) + \nabla(\zeta p) = \zeta f - \zeta u + p\, \nabla\zeta + \zeta \Delta u - \Delta(\zeta u) \in V_\beta^0(K)
\]
and $-\nabla\cdot (\zeta u) = \zeta g-u\cdot \nabla\zeta \in V_\beta^1(K)$.
By our assumption on $\beta$, the strip $-1/2\le \mbox{Re}\lambda \le -\beta+1/2$ does not contain eigenvalues
of the pencil ${\cal L}(\lambda)$. Using regularity results for solutions of the Stokes system
(see, e.g., \cite[Chapter 3, Theorem 5.6]{np-94} and \cite[Theorems 10.3.1 and 10.3.4]{mr-10}), we conclude that
$\zeta u \in V_\beta^2(K)$ and $\zeta p \in V_\beta^1(K)$. Since $\zeta$ has compact support, this implies
$\zeta(u,p) \in E_\beta^2(K) \times V_\beta^1(K)$.

It remains to prove that $\eta(u,p) \in E_\beta^2(K) \times V_\beta^1(K)$. Obviously,
\begin{equation} \label{1l7}
s\eta u - \Delta(\eta u) + \nabla(\eta p) = F, \quad -\nabla\cdot (\eta u) =G,
\end{equation}
where
$F=\eta f + p\, \nabla\eta + \eta \Delta u - \Delta(\eta u) \in V_\beta^0(K)$ and
$G = \eta g - u\cdot \nabla \eta \in X_\beta^1(K)$.
Suppose first that $\frac 12 - \lambda_1^+ < \beta \le 0$ (this requires that $\lambda_1^+>1/2$).
Then it follows from our assumptions that $\eta u \in V_{\beta-1}^1(K)$ and $\eta p \in V_{\beta-1}^0(K)$.
Applying Lemma \ref{l4}, we get $\eta(u,p) \in E_\beta^2(K) \times V_\beta^1(K)$.

Suppose now that $\max(0,-\lambda_1^+ +1/2) < \beta < 1/2$. Since $\beta>0$ and $F$ and $G$ are zero near the vertex of the cone,
we have $F \in E_0^0(K)$ and $G \in X_0^1(K)$. Furthermore, $\eta u \in V_{-1}^1(K)$ and $\eta p \in V_{-1}^0(K)$.
Thus, Lemma \ref{l4} implies $\eta(u,p) \in E_0^2(K)\times V_0^1(K)$. On the other hand,
$F\in E_\beta^0(K)$ and $G\in X_\beta^1(K)$. Since the interval $-\beta-1/2 \le \lambda \le -1/2$ is free of eigenvalues of the
pencil ${\cal N}(\lambda)$, we conclude from Lemma \ref{l5b} that $\eta(u,p) \in E_\beta^2(K) \times V_\beta^1(K)$.
The proof is complete. \hfill $\Box$ \\

As a consequence of the last two lemmas and Theorem \ref{t2}, we obtain the following result.

\begin{Th} \label{t3}
Suppose that $|s|=1$, $\mbox{\em Re}\, s \ge 0$ and $\frac 12 -\lambda_1^+ <\beta < \frac 12$.
Then the operator {\em (\ref{1t2})} is an isomorphism.
\end{Th}

P r o o f.
Since $0 < \frac 12 - \beta < \lambda_1^+$ and $-1 < -\frac 12 - \beta < \lambda_1^+ -1 \le 0$, the line
$\mbox{Re}\, \lambda = \frac 12 - \beta$ is free of eigenvalues of the pencil ${\cal L}(\lambda)$ and the line
$\mbox{Re}\, \lambda = -\frac 12 - \beta$ is free of eigenvalues of the pencil ${\cal N}(\lambda)$.
Thus by Theorem \ref{t2}, the operator (\ref{1t2}) has closed range. Furthermore by Lemma \ref{l6}, the kernel
of this operator is trivial. From Theorem \ref{t1} and Lemma \ref{l7} we conclude that the range of the operator
(\ref{1t2}) contains the set
\[
\big( E_0^{-1}(K)\cap E_\beta^0(K)\big)\times \big( L_2(K)\cap \big( V_0^1(K)\big)^* \cap X_\beta^1(K)\big)
\]
Since the range is closed it follows that the operator is an isomorphism onto
$E_\beta^0(K)\times X_\beta^1(K)$. \hfill $\Box$ \\

\noindent {\bf 1) The case \boldmath $\beta > 1/2$} \\

We consider the adjoint operator
\[
E_{-\beta}^0(K)\times \big( X_\beta^1(K)\big)^* \ni (v,q) \to A_\beta^*(v,q)
  \in \big( E_\beta^2(K)\cap \stackrel{\circ}{E}\!{}_\beta^1(K)\big)^* \cap \big( V_\beta^1(K)\big)^*
\]
of $A_\beta$,
where
\[
\langle A_\beta^*(v,q) \, ,\, (u,p)\rangle = \int_K\big( v\, (su-\Delta u+\nabla p)-q\, \nabla\cdot u\big)\, dx
\]
for all $u\in E_\beta^2(K)\cap \stackrel{\circ}{E}\!{}_\beta^1(K)$, $p\in V_\beta^1(K)$. We show that the
constant vector-function $(v,q)=(0,1)$ is an element of the kernel of $A_\beta^*$ if $1/2 < \beta < 5/2$.

\begin{Le}  \label{l11a}
The constant function $q=1$ belongs to $\big( X_\beta^1(K)\big)^*$ if and only if $1/2 < \beta < 5/2$. In this case, the pair
$(v,q)=(0,1)$ is an element of $\mbox{\em ker}\, A_\beta^*$, and $A_\beta$ is a mapping from
$\big( E_\beta^2(K)\cap \stackrel{\circ}{E}\!{}_\beta^1(K)\big)\times V_\beta^1(K)$ into the space
\begin{equation} \label{1l11a}
\Big\{ (f,g) \in E_\beta^0(K) \times X_\beta^1(K): \ \int_K g(x)\, dx =0\Big\}.
\end{equation}
\end{Le}

P r o o f.
Let $\zeta=\zeta(r)$ be a continuously differentiable function, $\zeta(r)=1$ for $r<1/4$, $\zeta(r)=0$ for $r>1/2$.
In the case $1/2 < \beta < 5/2$, the $V_{1-\beta}^0(K)$-norm of $\zeta$ and the
$V_{-\beta}^1(K)$-norm of $\eta=1-\zeta$ are finite, and we obtain
\[
\Big| \int_K g\, dx \Big| \le \Big| \int_K \zeta\, g\, dx \Big| + \Big| \int_K \eta \, g\, dx \Big|
  \le \| \zeta\|_{V_{1-\beta}^0(K)}\ \| g\|_{V_{\beta}^1(K)} + \| \eta\|_{V_{-\beta}^1(K)} \ \| g \|_{(V_{-\beta}^1(K))^*}
\]
for arbitrary $g\in X_\beta^1(K)$. In particular, this means that the constant function $q=1$ belongs to the dual space of $X_\beta^1(K)$
if $1/2 < \beta < 5/2$. Therefore, we get
\[
\int_K g\, dx = -\int_K q\,\nabla\cdot u\, dx =  \int_K u\, \nabla q\, dx =0
\]
for $q=1$, $u\in E_\beta^2(K) \cap\stackrel{\circ}{E}\!{}_\beta^1(K)$ and $g=-\nabla \cdot u$. In particular, $A_\beta^*(0,1)=0$.

We show that $1\not\in \big( X_\beta^1(K)\big)^*$ if $\beta\le 1/2$ or $\beta\ge 5/2$. Suppose that $\beta\ge 5/2$.
We consider the function $g(x)=\zeta(r)\, r^{-3}\, (\log r)^{-1}$. One easily checks that
$g\in V_\beta^1(K) \cap V_{\beta+1}^0(K) \subset X_\beta^1(K)$. However,
\[
\Big| \int_K g(x)\, dx \Big| \ge c\, \Big| \int_0^{1/4} \frac{dr}{r\log r} \Big| = \infty.
\]
In the case $\beta\le 1/2$, the function $g(x)=\chi(r)\, r^{-3}\, (\log r)^{-1}$, where $\chi(r)=0$ for $r<2$
and $\chi(r)=1$ for $r>4$, is an element of the space $V_\beta^1(K) \cap V_{\beta+1}^0(K)$, whereas the
integral of $g$ over $K$ is not finite. This shows that in both cases $1\not\in\big( X_\beta^1(K)\big)^*$.
The proof of the lemma is complete. \hfill $\Box$ \\

In order to prove the existence of solutions of the problem (\ref{par3}) in the space
$E_\beta^2(K)\times V_\beta^1(K)$ for $\beta>1/2$, we will apply the following regularity assertion.

\begin{Le} \label{l11}
Suppose that $(u,p) \in  E_\beta^2(K) \times V_\beta^1(K)$
is a solution of the problem {\em (\ref{par3})}, where
\[
f\in E_\beta^0(K) \cap E_\gamma^0(K), \quad g\in X_\beta^1(K)\cap X_\gamma^1(K),
\]
$-1/2 < \beta < \gamma < \mu_2^+ +1/2$, $\beta\not=1/2$ and $\gamma\not=1/2$.
In the case $\beta< 1/2 < \gamma$, we assume in addition that $\int_K g(x)\, dx =0$. Then $u\in E_\gamma^2(K)$ and $p\in V_\gamma^1(K)$.
\end{Le}

P r o o f.
In the cases $-1/2 < \beta < \gamma <1/2$ and $1/2 < \beta<\gamma < \mu_2^++1/2$,
the interval $-\gamma-1/2 \le \lambda \le -\beta-1/2$ does not contain eigenvalues of the pencil
${\cal N}(\lambda)$, and the assertion follows from Lemma \ref{l5b}. Therefore, we may assume that
$-1/2 < \beta < 1/2 < \gamma < \mu_2^++1/2$. In this case, the number $\lambda=-1$ is the
only eigenvalue of the pencil ${\cal N}(\lambda)$ in the interval $-\gamma-1/2 \le \lambda \le -\beta-1/2$.
Suppose first that $\gamma\le \beta+1/2$. Then $p$ is a solution of the problem (\ref{Neumann}),
where $F\in \big(V_{-\beta}^1(K)\big)^* \cap \big( V_{1-\gamma}^2(K)\big)^*$.
More precisely, $F$ has the representation (\ref{1l5b}) on $V_{-\beta}^1(K)$ and the representation (\ref{2l5b})
on $V_{1-\gamma}^2(K)$. As was shown in the proof of Lemma \ref{l5}, there exists a
pair $(v,\phi) \in V_{\gamma-1}^0(K) \times V_{\gamma-1}^{-1/2}(\partial K)$ such that
\begin{equation} \label{2l11}
- \int_K v\, \Delta q\, dx + \int_{\partial K} \phi\, \frac{\partial q}{\partial n}\, d\sigma = \langle F,q\rangle
  \ \mbox{ for all }q\in V_{1-\gamma}^2(K).
\end{equation}
By Lemma \ref{l5}, the function $p$ admits the decomposition
\[
p = c\, r^{-1} + v.
\]
Let $\zeta$ be the same cut-off function as in the proof of Lemma \ref{l5}, and let $\eta=1-\zeta$. Obviously,
$\zeta \in V_{1-\beta}^2(K)\subset V_{-\beta}^1(K)$ and $\eta \in V_{1-\gamma}^2(K)$. Furthermore,
$\frac{\partial \zeta}{\partial n} = \frac{\partial \eta}{\partial n}=0$ on $\partial K\backslash \{ 0\}$.
By (\ref{Neumann}) and Lemma \ref{l5a}, we have
\[
- \int_K p\, \Delta \zeta\, dx = \langle F,\zeta \rangle = \int_K \big( f-su-\nabla g\big)\cdot \nabla \zeta\, dx
  + \int_{\partial K} \sum_{i,j=1}^n \frac{\partial u_i}{\partial x_j}\, \Big( n_j\, \frac{\partial \zeta}{\partial x_i}
  - n_i\, \frac{\partial \zeta}{\partial x_j}  \Big)\, d\sigma,
\]
Furthermore, (\ref{2l11}) yields
\begin{eqnarray*}
&& \int_K \big( p-c\, r^{-1}\big)\, \Delta\zeta\, dx = - \int_K v\, \Delta\eta\, dx = \langle F, \eta \rangle \\
&& \quad = \ \int_K \Big( (f-\nabla g)\cdot \nabla \eta-sg\, \eta\Big)\, dx + \int_{\partial K} \sum_{i,j=1}^n
  \frac{\partial u_i}{\partial x_j}\, \Big( n_j\, \frac{\partial \eta}{\partial x_i} - n_i\, \frac{\partial \eta}{\partial x_j}
  \Big)\, d\sigma.
\end{eqnarray*}
Adding the last two equalities, we obtain
\[
- c\, \int_K r^{-1}\, \Delta \zeta\, dx = -s\, \int_K \big( \eta\, g +u\, \nabla\zeta\big)\, dx.
\]
Since $u \in \stackrel{\circ}{E}\!{}_\beta^1(K)$, $\nabla\cdot u = -g \in \big(V_{-\beta}^1(K)\big)^*$ and
$\zeta\in V_{-\beta}^1(K)$, this implies
\[
c\, \int_K r^{-1}\, \Delta \zeta\, dx = s\, \int_K g\, dx =0.
\]
The left-hand side of the last equality is equal to
\[
c\, \mbox{mes}\, \Omega \, \int_0^\infty \big( r\, \zeta''(r)+ 2\zeta'(r)\big)\, dr =
  - c\, \mbox{mes}\, \Omega,
\]
Hence, we get $c=0$ and $p= v\in V_{\gamma-1}^0(K)$. Furthermore, $u\in E_\beta^2(K) \subset V_{\gamma-1}^1(K)$.
Applying Lemma \ref{l4}, we obtain $u\in E_\gamma^2(K)$ and $p\in V_\gamma^1(K)$. This proves the lemma
for the case $\gamma\le \beta+1/2$. By induction, the assertion of the lemma holds for arbitrary
$\gamma \in (1/2\, ,\mu_2^+ +1/2)$. \hfill $\Box$ \\

Now, the following result holds analogously to Theorem \ref{t3}.

\begin{Th} \label{t4}
Suppose that $|s| =1$, $\mbox{\em Re}\, s \ge 0$, $\frac 12 <\beta < \min\big( \mu_2^+ +\frac 12\, , \, \lambda_1^+ +\frac 32\big)$.
Then the operator $(u,p) \to (f,g)$ of the problem {\em (\ref{par3})} is an isomorphism from
$\big( E_\beta^2(K) \cap\stackrel{\circ}{E}\!{}_\beta^1(K)\big)\times V_\beta^1(K)$ onto the space {\em (\ref{1l11a})}.
\end{Th}

P r o o f.
By Theorem \ref{t2} and Lemma \ref{l6}, the operator (\ref{1t2}) has closed range and trivial kernel.
Furthermore, the problem (\ref{par3}) is solvable in $E_\beta^2(K)\times V_\beta^1(K)$
for arbitrary $f\in C_0^\infty(\overline{K}\backslash \{ 0\})$ and $g\in C_0^\infty(\overline{K}\backslash \{ 0\})$
such that $\int_K g\, dx =0$. Indeed, for such $f$ and $g$, there exists a unique variational solution
$(u,p) \in  \stackrel{\circ}{E}\!{}_0^1(K) \times \big( L_2(K)+V_0^1(K)\big)$ of the problem (\ref{par3}).
Using Lemmas \ref{l7} and \ref{l11}, we conclude that $(u,p)\in E_\beta^2(K)\times V_\beta^1(K)$.
Since the set $C_0^\infty(\overline{K}\backslash \{ 0\})$ is dense in $E_\beta^0(K)$ and $X_\beta^1(K)$,
it follows that the problem (\ref{par3}) is solvable for arbitrary $f\in E_\beta^0(K)$, $g\in X_\beta^1(K)$,
$\int_K g\, dx=0$. \hfill $\Box$

\subsection{Necessity of the conditions on \boldmath $\beta$}

We already proved (see Lemmas \ref{l6b} and \ref{l6a}) that the operator $A_\beta$ is not Fredholm
if the line $\mbox{Re}\, \lambda = -\beta+1/2$ contains eigenvalues of the pencil ${\cal L}(\lambda)$ or
if $-\beta-1/2$ is an eigenvalue of the pencil ${\cal N}(\lambda)$. Our goal is to show that the conditions
$\beta>-\lambda_1^+ +1/2$ and $\beta<\min\big( \mu_2^+ +\frac 12\, , \, \lambda_1^+ +\frac 32\big)$ in Theorems
\ref{t3} and \ref{t4} are necessary.

\begin{Le} \label{l12}
Suppose that $-1/2 < \beta < -\lambda_1^+ +1/2$. Then the kernel of the adjoint operator $A_\beta^*$ is not trivial.
\end{Le}

P r o o f.
We assume that $-1/2 < \beta < -\lambda_1^+ +1/2$ and $\mbox{ker}\, A_\beta^*$ is trivial. Then it follows from Theorem \ref{t2}
and Lemma \ref{l6} that $A_\beta$ is an isomorphism. Let $(\phi,\psi)$
be an eigenvector of the pencil ${\cal L}(\lambda)$ corresponding to the eigenvalue $\lambda_1^+$ and let
$\zeta$ be a smooth ($\in C^2$) function with compact support equal to one near the vertex of the cone.
Furthermore, let $\delta$ be a real number, $-\lambda_1^+ +1/2 < \delta <\min(\beta+1,\, 1/2)$.
We consider the functions
\[
u(x)= \zeta(x)\, r^{\lambda_1^+}\, \phi(\omega), \quad p(x) = \zeta(x)\, r^{\lambda_1^+ -1}\, \psi(\omega).
\]
Obviously, $u\in E_\delta^2(K)$, $p\in V_\delta^1(K)$ and $u=0$ on $\partial K$.
Since
\[
-\Delta \, r^{\lambda_1^+}\, \phi + \nabla r^{\lambda_1^+ -1}\, \psi=0, \quad
-\nabla \cdot r^{\lambda_1^+}\, \phi =0 \ \mbox{ in }K,
\]
we get $su-\Delta u + \nabla p =f\in E_\beta^0(K)\cap E_\delta^0(K)$ and
$-\nabla\cdot u =g \in X_\beta^1(K)\cap X_\delta^1(K)$. Since $A_\beta$ is an isomorphism, there exists
a pair $(v,q) \in E_\beta^2(K)\times V_\beta^1(K)$ such that $A_\beta(v,q)=(f,q)$. From Lemma \ref{l11}
it follows that $(v,q)\in E_\delta^2(K) \times V_\delta^1(K)$. Therefore, $(u-v,p-q)\in \mbox{ker}\, A_\delta$.
Using Theorem \ref{t2}, we conclude that $(u,p)=(v,q)$, i.~e., $u\in E_\beta^2(K)$ and $p\in V_\beta^1(K)$.
Since obviously $(u,p) \not\in E_\beta^2(K) \times V_\beta^1(K)$, we got a contradiction. \hfill $\Box$ \\

Note that, in contrast to the case $1/2 < \beta < 5/2$, the kernel of $A_\beta^*$ does not contain the constant vector $(0,1)$
if $\beta$ satisfies the condition of Lemma \ref{l12}.

Next, we consider the case that $\beta > \min\big( \mu_2^+ +1/2\, , \, \lambda_1^+ +3/2\big)$.

\begin{Le} \label{l12a}
Suppose that $\lambda_1^+ < \mu_2^+ -1$. Then

{\em 1)} the kernel of the operator $A_\beta$ is not trivial for $\lambda_1^+ +3/2 < \beta < \mu_2^+ +1/2$.

{\em 2)} the kernel of the operator $A_\beta^*$ is not trivial for $-\mu_2^+ -1/2 < \beta < -\lambda_1^+ +1/2$.
\end{Le}

P r o o f.
1) Let $\big( \phi(\omega),\psi(\omega)\big)$ be an eigenvector to the eigenvalue $\lambda_1^-=-1-\lambda_1^+$ of the pencil
${\cal L}(\lambda)$. Furthermore, let $\zeta=\zeta(r)$ be a smooth function on $(0,\infty)$ with compact support which is equal
to one near the point $r=0$. We define $u_0(x)=\zeta(r)\, r^{\lambda_1^-}\phi(\omega)$ and
$p_0(x)=\zeta(r)\, r^{\lambda_1^- -1}\psi(\omega)$. One can easily check that $u_0 \in E_\beta^2(K)$ and $p_0 \in V_\beta^1(K)$
for $\beta>\lambda_1^+ +3/2$. Since
\[
-\Delta r^{\lambda_1^-}\phi(\omega) + \nabla r^{\lambda_1^- -1}\psi(\omega) =0 \ \mbox{ and } \
-\nabla\cdot r^{\lambda_1^-}\phi(\omega) =0 \ \mbox{ in } K
\]
we get $\big| su_0-\Delta u_0 + \nabla p_0\big| \le c\, r^{\lambda_1^-}$. Furthermore,  $\nabla\cdot u_0(x)=0$ for small
and for large $|x|$. Consequently, $s\, u_0 -\nabla u_0 + \nabla p_0 \in E_\gamma^0(K)$ and $\nabla\cdot u_0 \in X_\gamma^1(K)$
for arbitrary $\gamma \in (\frac 12 \, ,\, \lambda_1^+ + \frac 32)$. By Theorem \ref{t4}, there exists a vector function
$(u_1,p_1) \in E_\gamma^2(K) \times V_\gamma^1(K)$ such that
\[
s\, u_1 -\Delta u_1 + \nabla p_1=s\, u_0 -\Delta u_0 + \nabla p_0, \quad \nabla\cdot u_1=\nabla\cdot u_0\ \mbox{ in }K
\]
and $u_1=u_0=0$ on $\partial K$. Using Lemma \ref{l5b}, we conclude that $(u_1,p_1)\in E_\beta^2(K) \times V_\beta^1(K)$.
This means that $(u,p)=(u_0-u_1,p_0-p_1)$ is an element of the kernel of $A_\beta$. Since $(u_0,p_0) \not\in E_\gamma^2(K)
\times V_\gamma^1(K)$, the pair
 $(u,p)$ is not zero. This proves the first assertion.

2) We show that
\begin{equation} \label{1l12a}
\mbox{ker}\,A_\beta^* \supset \mbox{ker}\,A_{-\gamma} \quad\mbox{if } \gamma \le \beta\le \gamma+2.
\end{equation}
Let $(u,p)\in \big( E_{-\gamma}^2(K)\cap \stackrel{\circ}{E}\!{}_{-\gamma}^1(K)\big)\times V_{-\gamma}^1(K)$ be an element of the
kernel of $A_{-\gamma}$.  Then
\begin{equation} \label{2l12a}
\int_K \big( (su-\Delta u + \nabla p)\cdot v -q\, \nabla\cdot u\big)\, dx = 0 \ \mbox{ for all } v\in E_\gamma^0(K), \
  q\in \big( X_{-\gamma}^1(K)\big)^*.
\end{equation}
If $\gamma\le \beta \le \gamma+2$, then
\[
E_\beta^2(K) \subset E_\gamma^0(K) \quad\mbox{and}\quad E_{-\gamma}^2(K) \subset E_{-\beta}^0(K).
\]
Furthermore,
\[
V_\beta^1(K) \subset V_\gamma^1(K)+V_{\gamma+2}^1(K) \subset V_\gamma^1(K)+V_{\gamma+1}^0(K)
 \subset V_\gamma^1(K)+\big( V_{-\gamma}^1(K)\big)^* = \big( X_{-\gamma}^1(K)\big)^*
\]
and $V_{-\gamma}^1(K) \subset \big( X_\beta^1(K)\big)^*$. Thus, $(u,p) \in E_{-\beta}^0(K) \times \big( X_\beta^1(K)\big)^*$ and
the equality (\ref{2l12a}) is valid for all $v\in E_{\beta}^2(K)\cap \stackrel{\circ}{E}\!{}_{\beta}^1(K)$ and $q\in V_{\beta}^1(K)$.
Integrating by parts in (\ref{2l12a}), we get
\[
\int_K \big( u\, (sv-\Delta v+\nabla q)-p\, \nabla\cdot v\big)\, dx = 0 \quad\mbox{for all }
  v\in E_{\beta}^2(K)\cap \stackrel{\circ}{E}\!{}_{\beta}^1(K), \ q\in V_{\beta}^1(K)
\]
and, consequently, $(u,p)\in\mbox{ker}\,A_\beta^*$. This proves (\ref{1l12a}).

By the first assertion, the kernel of $A_{-\gamma}$ is not trivial for $-\mu_2^+ -1/2 < \gamma < -\lambda_1^+ -3/2$.
Consequently, by (\ref{1l12a}), the kernel of $A_\beta^*$ is not trivial for $-\mu_2^+ -1/2 < \beta < -\lambda_1^+ +1/2$.
\hfill $\Box$ \\

It remains to consider the case  $\lambda_1^+ > \mu_2^+ -1$. We construct some special solutions
of the system $(s-\Delta)\, u+\nabla p=0$, $\nabla\cdot u=0$. Let $\nu(x)$ denote the distance of the point
$x$ from the boundary $\partial K$. Obviously, the function $\nu$ is positively homogeneous of degree 1.
In the neighborhood $\nu(x)<\delta|x|$ of the boundary $\partial K$ with sufficiently small $\delta$,
the function $\nu$ is two times continuously differentiable and satisfies
the equality $|\nabla \nu|=1$. Furthermore, the vector $\nabla \nu(x)$ is orthogonal
to $\partial K$ in any point $x\in \partial K$.
For an arbitrary vector function $u$ in the neighborhood $\nu(x)<\delta |x|$ of $\partial K$, we define
\[
u_\nu = u\cdot \nabla\nu \quad\mbox{and}\quad u_\tau = u-u_\nu\, \nabla \nu.
\]
Obviously $u_\tau\cdot\nabla\nu=0$.

\begin{Le} \label{l12c}
Let $f$ and $g$ be positively homogeneous functions of degree $\mu$ in the neighborhood $\nu<\delta r$ of the boundary
$\partial K$. Furthermore, let $p$ and $q$ be polynomials of degree $k$ and $k+1$, respectively, such that
\[
p'(\nu)-\sqrt{s}\, p(\nu)=\nu^k, \quad q''(\nu)-2\sqrt{s}\, q'(\nu) = -\nu^k, \ q(0)=0.
\]
Then the functions
\[
U= e^{-\nu\sqrt{s}} \big( p(\nu)\,  g\, \nabla\nu + q(\nu)\, f_\tau\big), \quad P=e^{-\nu\sqrt{s}} p(\nu)\, f_\nu,
\]
satisfy the equalities
\[
(s-\Delta)U+\nabla P = e^{-\nu\sqrt{s}} \, \Big( \nu^k f + \big( p''(\nu)-2k\nu^{k-1}\big)\, g\nabla \nu + R\Big)
\]
and
\[
\nabla\cdot U = e^{-\nu\sqrt{s}} \, \big( \nu^k g + p(\nu)\, \, \nabla\cdot (g\, \nabla\nu) + q(\nu)\nabla\cdot f_\tau\big),
\]
in the neighborhood $\nu<\delta r$ of $\partial K$. Here $R$ is a finite sum of terms of the form
$\nu^j \phi$ with integer $j$, $0\le j\le k+1$, and positively homogeneous functions $\phi$ of degree $\mu-1$
or $\mu-2$.
\end{Le}

P r o o f.
Since $\nabla\, v(\nu)=v'(\nu)\, \nabla\nu$, we obtain
\begin{eqnarray*}
\nabla\cdot U & = & -\sqrt{s}\, e^{-\nu\sqrt{s}} \, p(\nu)\, g
  + e^{-\nu\sqrt{s}}\, \nabla \cdot \big( p(\nu)\,  g\, \nabla\nu + q(\nu)\, f_\tau\big) \\
& = & e^{-\nu\sqrt{s}} \, \Big( \nu^k \, g + p(\nu)\, \nabla\cdot\big(g\, \nabla\nu\big)
  + q(\nu)\nabla\cdot f_\tau\Big).
\end{eqnarray*}
Note that the functions $\nabla\cdot\big(g\, \nabla\nu\big)$ and $\nabla\cdot f_\tau$ are homogenous of degree $\mu-1$.
Furthermore,
\[
\nabla P =   f_\nu\, \nabla\big(  e^{-\nu\sqrt{s}}\, p(\nu)\big) + e^{-\nu\sqrt{s}}\, p(\nu)\, \nabla f_\nu
  = e^{-\nu\sqrt{s}} \, \big( \nu^k\, f_\nu\nabla\nu + p(\nu)\, \nabla f_\nu\big).
\]
Using the equalities $(s-\Delta)\, e^{-\nu\sqrt{s}} = \sqrt{s}\, e^{-\nu\sqrt{s}}\, \Delta\nu$ and
$\Delta p(\nu)=p''(\nu)+p'(\nu)\, \Delta \nu$, we get
\begin{eqnarray*}
(s-\Delta)\, U & = & \sqrt{s}\, e^{-\nu\sqrt{s}}\, \big( p(\nu)\, g\nabla\nu+ q(\nu)\, f_\tau\big)\, \Delta\nu \\
&&   + 2\sqrt{s}\, e^{-\nu\sqrt{s}} \sum_{j=1}^3 \frac{\partial \nu}{\partial x_j}\,
  \frac{\partial}{\partial x_j} \big( p(\nu)\, g\nabla\nu+ q(\nu)\, f_\tau\big)
 -e^{-\nu\sqrt{s}}\, \Delta\big( p(\nu)\, g\nabla\nu+ q(\nu)\, f_\tau\big)\\
& = & e^{-\nu\sqrt{s}}\, \Big( \nu^k\, f_\tau + (2\sqrt{s}\, p'(\nu)-p''(\nu)\big)\, g\nabla u +R\Big),
\end{eqnarray*}
where
\begin{eqnarray*}
R & = & p(\nu)\, \Delta(g\nabla\nu) + \big( \sqrt{s}\, p(\nu)-p'(\nu)\big)\,
  \Big( g\, \Delta\nu\, \nabla\nu + 2\sum_{j=1}^3 \frac{\partial\nu}{\partial x_j} \, \frac{\partial(g\nabla\nu)}{\partial x_j}\Big)\\
&&  + q(\nu)\, \Delta f_\tau + \big( \sqrt{s}\, q(\nu)-q'(\nu)\big)\,
  \Big( f_\tau\, \Delta\nu + 2\sum_{j=1}^3 \frac{\partial\nu}{\partial x_j} \,\frac{\partial f_\tau}{\partial x_j}\Big).
\end{eqnarray*}
Obviously, $\Delta(g\nabla\nu)$ and $\Delta f_\tau$ are positively homogeneous of degree $\mu-2$, while
$g\, \Delta\nu\, \nabla\nu$, $f_\tau\, \Delta\nu$, $\frac{\partial\nu}{\partial x_j} \, \frac{\partial(g\nabla\nu)}{\partial x_j}$
and $\frac{\partial\nu}{\partial x_j} \,\frac{\partial f_\tau}{\partial x_j}$ are positively homogeneous of degree $\mu-1$.
This proves the lemma. \hfill $\Box$

\begin{Le} \label{l12b}
Suppose that $\lambda_1^+ > \mu_2^+-1$. Then $\mbox{\em dim\,ker}\, A_\beta^* \ge 2$ for $\mu_2^+ +1/2 < \beta < \lambda_1^+ + 3/2$.
\end{Le}

P r o o f.
By Lemma \ref{l11a}, the pair $(v,q)=(0,1)$ belongs to the kernel of $A_\beta^*$. We construct another
element of $\mbox{ker}\, A_\beta^*$. Let $\phi$ be an eigenfunction of the pencil ${\cal N}(\lambda)$
corresponding to the eigenvalue $\lambda =\mu_2^+$. We start with the same pair $(u_0,p_0)$ as in the proof
of Lemma \ref{l6a}, i.e., $p_0=r^{\mu_2^+}\phi(\omega)$ and
\[
u_0(x) = v_0(x) - \chi\Big( \frac{\nu}r \Big) \, e^{-\nu\sqrt{s}}\, v_{0,\tau}(x),
\]
where $v_0= -s^{-1}\, \nabla p_0$, $v_{0,\tau} = v_0-v_{0,\nu}\nabla\nu$ for $\nu<\delta r$ and
$v_{0,\nu} = v_0\cdot \nabla\nu$ for $\nu<\delta r$. Here, $\chi$ is again a smooth cut-off function,
$\chi=1$ in $(0,\delta/2)$, $\chi=0$ in $(\delta,\infty)$. Then $u_0=0$ on $\partial K$ and
\[
(s-\Delta) u_0 + \nabla p_0  =  \chi\Big( \frac{\nu}r \Big) \, e^{-\nu\sqrt{s}} \, f_0   + F_0,  \quad
\nabla\cdot u_0 = \chi\Big( \frac{\nu}r \Big) \, e^{-\nu\sqrt{s}} \, g_0   + G_0
\]
in $K$, where
\[
f_0=-\sqrt{s}\, v_{0,\tau}\Delta\nu- 2\sqrt{s}\sum_{j=1}^3 \frac{\partial\nu}{\partial x_j}\, \frac{\partial v_{0,\tau}}{\partial x_j}\, ,
  \quad g_0 = -\nabla \cdot v_{0,\tau}
\]
and
\[
F_0 = \chi\Big( \frac{\nu}r \Big) \, e^{-\nu\sqrt{s}} \Delta v_{0,\tau}+ \Big[ \Delta,\chi\Big(\frac{\nu}r\Big)\Big]\,
  e^{-\nu\sqrt{s}}\, v_{0,\tau}, \quad G_0= - e^{-\nu\sqrt{s}}\, v_{0,\tau}\cdot \nabla \chi\Big( \frac{\nu}r \Big)
\]
(see (\ref{2l6a}) and (\ref{1l6a})). Obviously, the functions $f_0$ and $g_0$ are positively homogeneous of degree $\mu_2^+-2$.
Let $\eta$ be a smooth function on $\overline{K}$ such that $\eta(x)=0$ for $|x|<1$
and $\eta(x)=1$ for $|x|>2$. Then
\[
\eta u_0 \in E_{-\beta}^2(K), \quad \eta p_0 \in V_{-\beta}^1(K)
\]
since $\beta>\mu_2^++1/2$. Furthermore, $\eta F_0 \in V_\gamma^0(K)$ and $\eta G_0 \in X_\gamma^1(K)$ for arbitrary
$\gamma <\frac 32 -\mu_2^+$.

We define
\[
w_1 = \frac{1}{2\sqrt{s}}\, \chi\Big( \frac{\nu}r \Big) \, e^{-\nu\sqrt{s}} \, \big( 2g_0\nabla\nu-\nu\, f_{0,\tau}\big), \quad
q_1 = \frac{1}{\sqrt{s}}\, \chi\Big( \frac{\nu}r \Big) \, e^{-\nu\sqrt{s}} \, f_{0,\nu}\, .
\]
By Lemma \ref{l12c},
\[
(s-\Delta)\, w_1 + \nabla q_1 = - \chi\Big( \frac{\nu}r \Big) \, e^{-\nu\sqrt{s}} \, f_0 + F_1
\]
and
\[
\nabla\cdot w_1 = \chi\Big( \frac{\nu}r \Big) \, e^{-\nu\sqrt{s}} \, \big( -g_0 +g_1+\nu g_2\big) + G_1
\]
where $g_1 = s^{-1/2}\nabla\cdot(g_0\nabla\nu)$ and $g_2=-\frac 12 s^{-1/2}\nabla\cdot f_{0,\tau}$. Here,
$F_1$ is a sum of terms of the form
\begin{equation} \label{1l12b}
\chi\Big( \frac{\nu}r \Big) \, e^{-\nu\sqrt{s}} \, \Phi \quad\mbox{or}\quad
\chi^{(k)}\Big( \frac{\nu}r \Big) \, e^{-\nu\sqrt{s}} \, \Psi,
\end{equation}
where $\Phi$ is positively homogeneous of degree $\le \mu_2^+-3$ and $\Psi$ is positively homogeneous of degree $\le \mu_2^+-2$,
while $G_1$ is a sum of terms of the form (\ref{1l12b}), where  $\Phi$ is positively homogeneous of degree $\le \mu_2^+-4$.
Consequently, $\eta F_1 \in V_\gamma^0(K)$ and $\eta G_1 \in X_\gamma^1(K)$ for arbitrary $\gamma <\frac 32 -\mu_2^+$.
Obviously, $w_{1,\tau}=0$ and $w_{1,\nu}=s^{-1/2}g_0$ on $\partial K$. Since $g_0$ is positively homogeneous of
degree $\mu_2^+ -2$, there exists a solution of the problem
\[
\Delta p_1=0 \ \mbox{ in }K, \quad \frac{\partial p_1}{\partial n} = sw_{1,\nu}=s^{1/2}g_0 \ \mbox{ on }\partial K\backslash\{ 0\}
\]
which has the form $p_1=r^{\mu_2^+-1}\psi(\omega)$ if $\mu_2^+ -1$ is not an eigenvalue
or $p_1=r^{\mu_2^+-1}\big( \psi_1(\omega)+\psi_2(\omega)\log r\big)$ if $\mu_2^+ -1$ is an eigenvalue of the pencil
${\cal N}(\lambda)$ (see, e.~g., \cite[Lemma 6.1.13]{kmr-97}). We suppose for simplicity that $\mu_2^+ -1$ is not an eigenvalue
of the pencil ${\cal N}(\lambda)$ and set $v_1=-s^{-1}\nabla p_1$ and
\[
u_1(x) =v_1(x)- \chi\Big( \frac{\nu}r \Big) \, e^{-\nu\sqrt{s}} \, v_{1,\tau}\, .
\]
Then $(s-\Delta)u_1 + \nabla p_1 = F_2$ and
\[
\nabla\cdot u_1 = - \chi\Big( \frac{\nu}r \Big) \, e^{-\nu\sqrt{s}} \, \nabla\cdot v_{1,\tau} + G_2
\]
where $F_2$ is a sum of terms of the form (\ref{1l12b}) with positively homogeneous functions $\Phi$ and $\Psi$
of degree $\le \mu_2^+-3$ and degree $\le \mu_2^+-2$, respectively, and $G_2$ is a sum of terms of the form (\ref{1l12b})
with positively homogeneous functions $\Phi$ of degree $\le \mu_2^+-4$. In particular, $\eta F_2 \in V_\gamma^0(K)$
and $\eta G_2 \in X_\gamma^1(K)$ for arbitrary $\gamma <\frac 32 -\mu_2^+$. Obviously, $u_{1,\tau}=0$ and
$u_{1,\nu}=v_{1,\nu}=-s^{-1}\nabla p_1\cdot\nabla\nu=-w_{1,\nu}$ on $\partial K\backslash\{ 0\}$. This means that $u_1+w_1=0$
on $\partial K\backslash\{ 0\}$.

In the last step, we consider the function
\[
w_2 = s^{-1/2}\, \chi\Big( \frac{\nu}r \Big) \, e^{-\nu\sqrt{s}} \, \Big( \big( g_1-\nabla\cdot v_{1,\tau}\big)\, \nabla\nu
  + \big(\nu+s^{-1/2}\big)g_2\nabla\nu\Big).
\]
Using Lemma \ref{l12c}, we conclude that
\[
(s-\Delta)\, w_2 = F_3 \ \mbox{ and } \ \nabla\cdot w_2 = - \chi\Big( \frac{\nu}r \Big) \, e^{-\nu\sqrt{s}} \,
  \big( g_1-\nabla\cdot v_{1,\tau}\big) +\nu\, g_2\big) + G_3,
\]
where $\eta F_3 \in V_\gamma^0(K)$ and $\eta G_3 \in X_\gamma^1(K)$ for arbitrary $\gamma <\frac 32 -\mu_2^+$.
Furthermore, $w_{2,\tau}=0$ and $w_{2,\nu}=s^{-1/2}(g_1-\nabla\cdot v_{1,\tau}) + s^{-1}g_2$ on
$\partial K \backslash\{ 0\}$. Since the functions $g_1,g_2$ and $\nabla\cdot v_{1,\tau}$ are positively homogeneous
of degree $\mu_2^+-3$, there exists a solution $p_2$ of the problem
\[
\Delta p_2=0 \ \mbox{ in }K, \quad \frac{\partial p_2}{\partial n}=sw_{2,\nu} =  s^{1/2}(g_1-\nabla\cdot v_{1,\tau}) + g_2
   \ \mbox{ on }\partial K\backslash\{ 0\}
\]
which has the form $p_2=r^{\mu_2^+-2}\psi(\omega)$ if $\mu_2^+ -1$ or $p_2=r^{\mu_2^+-2}\big( \psi_1(\omega)+\psi_2(\omega)\log r\big)$.
We set $v_2=-s^{-1}\nabla p_2$ and
\[
u_2(x) =v_2(x)- \chi\Big( \frac{\nu}r \Big) \, e^{-\nu\sqrt{s}} \, v_{2,\tau}(x) .
\]
Then $\eta(s-\Delta)u_2+\nabla p_2 \in V_\gamma^0(K)$ and $\eta \nabla\cdot u_2 \in X_\gamma^1(K)$ for arbitrary
$\gamma <\frac 32 -\mu_2^+$. Furthermore $u_2=-w_2$ on $\partial K\backslash \{ 0\}$.

We define
\[
U = u_0+u_1+u_2+w_1+w_2 \quad\mbox{and}\quad P=p_0+p_1+p_2\, .
\]
Since $u_0$, $u_1+w_1$ and $u_2+w_2$ are zero on the boundary, the function $U$ also vanishes on
$\partial K\backslash\{ 0\}$. Furthermore, as was show above, we have
$\eta(s-\Delta)U + \eta P \in V_\gamma^0(K)$ and $\eta \nabla\cdot U \in X_\gamma^1(K)$ for arbitrary
$\gamma <\frac 32 -\mu_2^+$. Since $\nabla\eta=0$ outside the region $1<|x|<2$, we conclude from this that
\[
(s-\Delta)\, (\eta U) + \nabla(\eta P) \in V_\gamma^0(K), \quad \nabla\cdot (\eta U) \in X_\gamma^1(K)
\]
for arbitrary $\gamma <\frac 32 -\mu_2^+$. By the condition of the lemma on $\lambda_1^+$ and $\mu_2^+$,
the number $\gamma$ can be chosen such that
$\frac 12 -\lambda_1^+ < \gamma < \min(\frac 12,2-\beta)$. In particular, $\gamma>-\frac 12$.
By Theorem \ref{t3}, there exists a solution $(u_3,p_3) \in E_\gamma^2(K)\times V_\gamma^1(K)$
of the problem
\[
(s-\Delta)\, u_3+ \nabla p_3 = (s-\Delta)\, (\eta U) + \nabla (\eta P), \ \ \nabla\cdot u_3 = \nabla\cdot(\eta U)
\ \mbox{ in K}
\]
with the Dirichlet condition $u_3=0$ on $\partial K\backslash \{ 0\}$. This means that $u= \eta U -u_3$
and $p=\eta P -p_3$ satisfy the equations $(s-\Delta)u+\nabla p =0$, $\nabla\cdot u=0$ in $K$
and the boundary condition $u=0$ on $\partial K\backslash \{ 0\}$. Obviously, $u\not=0$ since
$\eta U \not\in E_\gamma^2(K)$. Furthermore,
\[
\eta(U,P)\in E_{-\beta}^2(K)\times V_{-\beta}^1(K) \subset E_{-\beta}^0(K) \times \big( X_\beta^1(K)\big)^*.
\]
The same is true for $(u_3,p_3)$. Indeed, we have $u_3 \in E_{\gamma-2}^0(K)\cap E_\gamma^0(K) \subset E_{-\beta}^0(K)$,
since $\gamma-2<-\beta<-\frac 12 < \gamma$, and from $p_3\in V_\gamma^1(K)$ it follows that
$\eta p_3 \in V_{-\beta}^1(K)$ (since $-\beta<\gamma$)
and $(1-\eta)p_3 \in V_{1-\beta}^0(K) \subset \big( V_\beta^1(K)\big)^*$ (since $1-\beta>\gamma-1$). Thus,
$p_3 \in V_{-\beta}^1(K)+\big( V_\beta^1(K)\big)^* = \big( X_\beta^1(K)\big)^*$.
This proves that $(u,p)\in E_{-\beta}^0(K) \times \big( X_\beta^1(K)\big)^*$ is an element of the kernel of $A_\beta^*$.
The proof is complete. \hfill $\Box$

\subsection{An estimate for arbitrary \boldmath $s\not=0, \ \mbox{Re}\, s \ge 0$}

We consider the problem (\ref{par3}), where $s$ is an arbitrary complex parameter $\not= 0$ with nonnegative real part.
As a consequence of Theorems \ref{t3} and \ref{t4}, we obtain the following result.

\begin{Th} \label{t5}
Suppose that $s \not=0$ and $\mbox{\em Re}\, s \ge 0$. Furthermore, we assume that $f\in V_\beta^0(K)$ and $g\in X_\beta^1(K)$, where
$ -\lambda_1^+ +1/2<\beta < \min\big( \mu_2^+ +1/2\, , \, \lambda_1^+ +3/2\big)$ and $\beta\not= 1/2$.
In the case $\beta>1/2$, we assume in addition that $\int_K g(x)\, dx =0$. Then there exists a uniquely determined
solution $(u,p) \in E_\beta^2(K) \times V_\beta^1(K)$ of the problem {\em (\ref{par3})}. This solution satisfies the estimate
\[
\| u\|_{V_\beta^2(K)} + |s|\, \| u\|_{V_\beta^0(K)} + \| p\|_{V_\beta^1(K)} \le c\, \Big( \| f\|_{V_\beta^0(K)}
  + \| g \|_{V_\beta^1(K)} + |s|\, \| g\|_{(V_{-\beta}^1(K))^*}\Big)
\]
with a constant $c$ independent of $f$, $g$ and $s$.
\end{Th}

P r o o f.
The pair $(u,p)$ is a solution of the problem (\ref{par3}) if and only if
\[
v(x) = u\big(|s|^{-1/2}x\big) \ \mbox{ and } \ q(x)=|s|^{-1/2}\, p\big(|s|^{-1/2}x\big)
\]
satisfy the equations
\begin{equation} \label{2t5}
\frac{s}{|s|}\, v - \Delta v + \nabla q = F, \ - \nabla\cdot v = G \ \mbox{ in }K, \quad v|_{\partial K\backslash \{ 0\}}=0,
\end{equation}
where $F(x)=|s|^{-1}\, f\big(|s|^{-1/2}x\big)$ and $G(x)=|s|^{-1/2}\, g\big(|s|^{-1/2}x\big)$.
By Theorems \ref{t3} and \ref{t4}, the boundary value
problem for the system (\ref{2t5}) has a uniquely determined solution $(v,q) \in E_\beta^2(K)\times V_\beta^1(K)$
satisfying the estimate
\[
\| v\|^2_{E_\beta^2(K)} + \| q\|^2_{V_\beta^1(K)} \le c\, \Big( \| F\|^2_{V_\beta^0(K)} + \| G\|^2_{V_\beta^1(K)}
  +  \| G\|^2_{(V_{-\beta}^1(K))^*}\Big).
\]
One easily verifies that
\[
 \| F\|^2_{V_\beta^0(K)} + \| G\|^2_{V_\beta^1(K)}
  +  \| G\|^2_{(V_{-\beta}^1(K))^*} = |s|^{\beta-1/2}\Big( \| f\|^2_{V_\beta^0(K)} + \| g\|^2_{V_\beta^1(K)}
  +  |s|^2\, \| g\|^2_{(V_{-\beta}^1(K))^*}\Big)
\]
and
\begin{eqnarray*}
\| v\|^2_{E_\beta^2(K)} + \| q\|^2_{V_\beta^1(K)} & \ge & \| v\|^2_{V_\beta^2(K)} + \| v\|^2_{V_\beta^0(K)} + \| q\|^2_{V_\beta^1(K)} \\
  & = & |s|^{\beta-1/2}\, \Big( \| u\|^2_{V_\beta^2(K)} + |s|^2\, \| u\|_{V_\beta^0(K)} + \| p\|^2_{V_\beta^1(K)} \Big).
\end{eqnarray*}
This proves the theorem. \hfill $\Box$

\section{The time-dependent problem in \boldmath $K$}

Now, we consider the problem (\ref{stokes1}), (\ref{stokes2}). Using the last theorem, we can easily show
that this problem has a uniquely determined solution in a certain weighted Sobolev space.

\subsection{Weighted Sobolev spaces in \boldmath $K\times{\Bbb R}_+$}

Let $Q=K\times {\Bbb R}_+ = K\times (0,\infty)$.
We denote by $W_\beta^{2l,l}(Q)$ the weighted Sobolev space of all functions $u=u(x,t)$ on $Q$ with finite norm
\[
\| u\|_{W_\beta^{2l,l}(Q)} = \Big( \int_0^\infty  \sum_{k=0}^l \| \partial_t^k u(\cdot,t) \|^2_{V_\beta^{2l-2k}(K)}\, dt \Big)^{1/2}.
\]
In particular, $W_\beta^{0,0}(Q) = L_2\big({\Bbb R}_+,V_\beta^0(K)\big)$ and
$W_\beta^{2,1}(Q)$ is the set of all $u\in L_2({\Bbb R}_+,V_\beta^2(K))$ such that $\partial_t u \in L_2({\Bbb R}_+,V_\beta^0(K))$.
The space $\stackrel{\circ}{W}\!{}_\beta^{2l,l}(Q)$ is the subspace of all $u\in W_\beta^{2l,l}(Q)$ satisfying
the condition $\partial_t^k u|_{t=0}$ for $x\in K$, $k=0,\ldots,l-1$. Note that
$\partial_t^k u(\cdot,0) \in V_\beta^{2l-2k-1}(K)$ for $u\in W_\beta^{2l,l}(Q)$,
$k=0,\ldots,l-1$ (see \cite[Proposition 3.1]{kozlov-89}).
By \cite[Proposition 3.4]{kozlov-89}, the Laplace transform realizes an isomorphism from $\stackrel{\circ}{W}\!{}_\beta^{2l,l}(Q)$
onto the space $H_\beta^{2l}$ of all holomorphic functions $\tilde{u}(x,s)$ for $\mbox{Re}\, s >0$ with values
in $E_\beta^{2l}(K)$ and finite norm
\[
\| \tilde{u}\|_{H_\beta^{2l}} = \sup_{\gamma>0} \Big( \int_{-\infty}^{+\infty}
  \sum_{k=0}^l |s|^{2k}\, \| \tilde{u}(\cdot,\gamma+i\tau)\|^2_{V_\beta^{2l-2k}(K)} \, d\tau  \Big)^{1/2}.
\]
The proof for the analogous result in nonweighted spaces can be found in \cite[Theorem 8.1]{av-64}.

\subsection{Existence and uniqueness of solutions}

As a consequence of Theorem \ref{t5}, we obtain the following assertion.

\begin{Th} \label{t6}
Suppose that $f\in L_2\big( {\Bbb R}_+,V_\beta^0(K)\big)$, $g\in L_2\big({\Bbb R}_+,V_\beta^1(K)\big)$
and $\partial_t g\in L_2\big({\Bbb R}_+,(V_{-\beta}^1(K))^*\big)$, where
$-\lambda_1^+ +1/2<\beta < \min\big( \mu_2^+ +1/2\, , \, \lambda_1^+ +3/2\big)$ and $\beta\not= 1/2$.
In the case $\beta>1/2$, we assume in addition that $\int_K g(x,t)\, dx =0$ for almost all $t$.
Then there exists a uniquely determined solution $(u,p) \in W_\beta^{2,1}(Q) \times L_2\big({\Bbb R}_+, V_\beta^1(K)\big)$
of the problem {\em (\ref{stokes1}), (\ref{stokes2})} satisfying the estimate
\[
\| u\|_{W_\beta^{2,1}(Q)} + \| p\|_{L_2({\Bbb R}_+,V_\beta^1(K)}  \le c\, \Big( \| f\|_{W_\beta^{0,0}(Q)}
  + \| g \|_{L_2({\Bbb R}_+,V_\beta^1(K))} +  \| \partial_t g\|_{L_2({\Bbb R}_+,(V_{-\beta}^1(K))^*)}\Big)
\]
with a constant $c$ independent of $f$, $g$.
\end{Th}

P r o o f.
Let $\tilde{f}\in H_\beta^0$ and $\tilde{g}$ be the Laplace transforms (with respect to the variable $t$) of $f$ and $g$.
For arbitrary $s\not=0$, $\mbox{Re}\, s\ge 0$, there exists a uniquely determined solution
$\tilde{u}(\cdot,s), \tilde{p}(\cdot,s) \in E_\beta^2(K)\times V_\beta^1(K)$ of the problem (\ref{par1})
satisfying the estimate
\begin{eqnarray*}
&& \| \tilde{u}(\cdot,s)\|^2_{V_\beta^2(K)} + |s|^2\, \| \tilde{u}(\cdot,s)\|^2_{V_\beta^0(K)}
  + \| \tilde{p}(\cdot,s)\|^2_{V_\beta^1(K)} \\
&& \le c\, \Big( \| \tilde{f}(\cdot,s)\|^2_{V_\beta^0(K)}
  + \| \tilde{g}(\cdot,s) \|^2_{V_\beta^1(K)} + |s|^2\, \| \tilde{g}(\cdot,s)\|^2_{(V_{-\beta}^1(K))^*}\Big)
\end{eqnarray*}
with a constant $c$ independent of $s$. Integrating over the line $\mbox{Re}\, s = \gamma$ and taking the supremum
with respect to $\gamma>0$, we get the assertion of the theorem. \hfill $\Box$ \\

Furthermore, we can deduce the following regularity result from Lemma \ref{l11}.

\begin{Th} \label{t7}
Suppose that $(u,p) \in  W_\beta^{2,1}(Q) \times L_2\big( {\Bbb R}_+,V_\beta^1(K)\big)$
is a solution of the problem {\em (\ref{stokes1}), (\ref{stokes2})}, where
\[
f\in L_2\big({\Bbb R}_+,V_\beta^0(K)\big) \cap L_2\big({\Bbb R}_+,V_\gamma^0(K)\big),
\]
\[
g\in L_2\big({\Bbb R}_+,V_\beta^1(K)\big) \cap L_2\big({\Bbb R}_+,V_\gamma^1(K)\big), \quad
\partial_t g\in L_2\big({\Bbb R}_+,(V_{-\beta}^1(K))^*\big) \cap L_2\big({\Bbb R}_+,(V_{-\gamma}^1(K))^*\big),
\]
$1/2-\lambda_1^+ < \beta,\gamma < \min(\mu_2^+ +1/2,\lambda_1^+ + 3/2)$, $\beta\not=1/2$ and $\gamma\not=1/2$.
In the case $\gamma>1/2$ we assume in addition that $\int_K g(x,t)\, dx =0$ for almost all $t$.
Then $u\in W_\gamma^{2,1}(Q)$ and $p\in L_2\big({\Bbb R}_+,V_\gamma^1(K)\big)$.
\end{Th}

P r o o f.
By Theorem \ref{t6}, there exists a uniquely determined solution $(v,q) \in W_\gamma^{2,1}(Q) \times
L_2\big({\Bbb R}_+,V_\gamma^1(K)\big)$ of the problem (\ref{stokes1}), (\ref{stokes2}). The Laplace transforms
$\big( \tilde{u}(x,s),\tilde{p}(x,s)\big)$ and $\big( \tilde{v}(x,s),\tilde{p}(v,s)\big)$ belong to the
spaces $E_\beta^2(K)\times V_\beta^1(K)$ and $E_\gamma^2(K)\times V_\gamma^1(K)$ for almost all $s$ in the
half-plane $\mbox{Re}\, s \ge 0$. Using Lemma \ref{l11}, we obtain $(\tilde{u},\tilde{p})
=(\tilde{v},\tilde{q})$ and, consequently, $(u,p)=(v,q)$. \hfill $\Box$

\end{document}